\documentclass[12pt,notitlepage,twoside,a4paper]{amsart}
 \usepackage{amsfonts}

\usepackage{amsmath,amssymb,enumerate}

\usepackage{epsfig,fancyhdr,color}

\usepackage{amssymb}
\usepackage{amsmath,amsthm}
\usepackage{latexsym}
\usepackage{amscd}
\usepackage{psfrag}
\usepackage{graphicx}
\usepackage[latin1]{inputenc}
\usepackage[all]{xy}
\usepackage[mathcal]{eucal}

\definecolor{NoteColor}{rgb}{1,0,0}

\renewcommand{\textsc}{\textcolor{red}}

\newtheorem{theorem}{\rm\bf Theorem}[section]
\newtheorem{proposition}[theorem]{\rm\bf Proposition}

\newtheorem*{theorem 1}{\rm\bf Proposition 1}
\newtheorem*{theorem 2}{\rm\bf Proposition 2}

\theoremstyle{definition}
\newtheorem{definition}[theorem]{\rm\bf Definition}

\theoremstyle{remark}
\newtheorem{remark}[theorem]{\rm\bf Remark}

\newtheorem{problem}[theorem]{\rm\bf Problem}

\def\interieur#1{\mathord{\mathop{\kern 0pt #1}\limits^\circ}}

\title[A Commentary on Teichm\"uller's paper]{A commentary on Teichm\"uller's paper \emph{Extremale quasikonforme Abbildungen und quadratische Differentiale}}

\author[Vincent Alberge, Athanase Papadopoulos and Weixu Su]{Vincent Alberge, Athanase Papadopoulos and Weixu Su}

\thanks{The authors were partially supported by the French ANR program FINSLER and by the U.S. National Foundation grants DMS 1107452, 1107263, 1107367 ``RNMS  GEometric structures And Representation varieties."}

 \address{
Institut de Recherche Math\'ematique Avanc\'ee\\ CNRS et Universit\'e de Strasbourg\\\small 7 rue Ren\'e
  Descartes - 67084 Strasbourg Cedex, France\\
  Department of Mathematics, Fudan University, 200433, Shanghai, China\\
Emails: alberge@math.unistra.fr, papadop@math.unistra.fr\\
  suwx@fudan.edu.cn}

\date{\today}

\begin{document}

  \maketitle
  \begin{abstract}

We provide a commentary on Teichm\"uller's paper
\emph{Extremale quasikonforme Abbildungen und quadratische Differentiale}
(Extremal  quasiconformal mappings of closed oriented Riemann surfaces),
Abh. Preuss. Akad. Wiss., Math.-Naturw. Kl. 1940, No.22, 1-197 (1940).  The paper is quoted in several works, although it was read by very few people. Some of the results it contains were rediscovered later on and published without any reference to Teichm\"uller.
In this commentary, we highlight the main results and the main ideas contained in that paper and we describe some of the important developments they gave rise to.

The final version of this paper, together with the English translation of Teichm\"uller's paper,  will apper in Volume V of the \emph{Handbook of Teichm\"uller theory} (European Mathematical Society Publishing House, 2015).

\end{abstract}

  \tableofcontents

%
%
%
%
%
%
%
%

\section{Introduction to Teichm\"uller's paper}

The paper \cite{T20}, published in 1939, is the most quoted article by Teichm\"uller although, like the rest of his papers, it was read by very few people. It contains the foundations of what we call now \textit{the classical Teichm\"uller theory}. The present notes are intended to be a sort of a \textit{reading guide} for that paper. In order to understand its logic, and especially to see what the author proves and what he admits without proof, the paper has to be read completely. Indeed, some of the results are first referred to as ``conjectures," but they are proved later on in the paper. Some of them are proved in subsequent papers. For example, the proof of the important \emph{Existence Theorem} for extremal quasiconformal mappings\index{extremal quasiconformal mappings!existence} is given in \cite{T29}, published in 1943. For the present commentary, we have used the English translation by G. Th\'eret which appears in the present volume.

We start with a list of major ideas contained in this paper.

\begin{enumerate}

\item \label{t101} A topology and a metric on Teichm\"uller and moduli spaces, with a formula for the dimensions of these spaces which coincides with Riemann's count of moduli. This is a solution to the so-called \emph{Riemann moduli problem},\index{Riemann moduli problem}\index{problem!Riemann moduli} which asks for a precise meaning of Riemann's statement that there are $3g-3$ ``moduli" for a closed Riemann surface of genus $g\geq 2$.

\item \label{t102} The idea of marking Riemann surfaces by mapping classes (homeomorphisms defined up to homotopy) from a fixed topological surface to a varying Riemann surface, in order to overcome problems caused by the singularities of moduli space.

\item \label{t000} The introduction of a cover of Riemann's moduli space which is nonsingular, on which the mapping class group acts properly discontinuously. This is the space called now \emph{Teichm\"uller space}. It is the space of equivalence classes of marked Riemann surfaces. The Riemann moduli space appears as a quotient of Teichm\"uller space by the action of the mapping class group.

\item  \label{t2} The use of quasiconformal mappings as a tool for understanding \emph{conformal} mappings and for the study of moduli of Riemann surfaces, and not only as a generalization of the notion of conformal mappings. In other words, in Teichm\"uller's works, quasi-conformality is not a notion of non-conformality but rather a tool for constructing a conformal structure for which this map is conformal.

\item \label{t218} The use of holomorphic quadratic differentials, their trajectory structure and the singular flat metric they induce on the surface as a new essential tool in the theory of moduli.

\item \label{t23} The introduction of the Teichm\"uller metric, the proof that this metric is Finsler and the study of its infinitesimal norm and its geodesics.

\item \label{t25} The study of the infinitesimal  theory of quasiconformal mappings, and a thorough study of the partial differential equations which they satisfy.

\item \label{t259} The translation of problems concerning conformal structures on surfaces into problems concerning Riemannian metrics, and the use of these metrics to solve problems in the theory of conformal mappings.

\item \label{t206}
 The introduction of an equivalence relation on the set of Beltrami differentials  and the identification of the space of equivalence classes as the tangent space to Teichm\"uller space at a point. This contains at the same time the idea of an almost-complex structure on the tangent space to Teichm\"uller space and it is at the basis of the theory of variation of higher-dimensional complex manifolds developed by Kodaira and Spencer.

\item \label{t208}  The identification of the tangent space at a point of Teichm\"uller space with the topological dual of the space of holomorphic quadratic differentials.

\item \label{t24} The introduction and the study of Teichm\"uller discs (which Teichm\"uller calls ``complex geodesics").

\item \label{t26} The comparison between the length of closed geodesics and the quasiconformal dilatation of a map between hyperbolic surfaces (rediscovered by Sorvali and by Wolpert).

\item \label{t27} A precise use of the method of continuity in the setting of the  moduli problem, into which Poincar\'e, Klein and others had bumped over a long period of time. In Teichm\"uller's setting, the method is applied (as it should be) to maps between objects which are known to be manifolds of the same dimension.

\item \label{t215} The question of whether there is a Hermitian metric on Teichm\"uller space.

\item \label{t218} A study of convexity properties of the Teichm\"uller metric and the question of studying totally geodesic subspaces for that metric.

 \item \label{t218}  The translation of Nielsen's realization problem into a question of finding a fixed point of the action of a finite subgroup of the mapping class group on Teichm\"uller space, and the appeal to the convexity properties for the solution. This is the approach that Kerckhoff used in the solution of the problem \cite{kerckhoff1}.
 
\item The introduction of the theory of \emph{non-reduced Teichm\"uller spaces},\index{non-reduced Teichm\"uller space}\index{Teichm\"uller space!non-reduced} as a theory where each point on the boundary of a surface (or on a union of arcs on this boundary), is considered as a distinguished point. In particular, he considers the non-reduced Teichm\"uller space of the disc (the universal Teichm\"uller space). 

\item The idea of a new metric defined on a Riemann surface, where the distance between two points is the logarithm of the least quasiconformal constant (or the ``dilatation quotient") of a self-map of the surface sending one point to the other.

\item The study of boundary maps of quasiconformal homeomorphisms of the disc.
\end{enumerate}

  For a more global view of Teichm\"uller's contribution on moduli of Riemann surfaces, the reader is also referred to his other papers \cite{T200},  \cite{T23},  \cite{T24}, \cite{T29}, \cite{T31}, \cite{T32} and the corresponding commentaries \cite{T200C}, \cite{T23C}, \cite{T24C}, \cite{T29C}, \cite{T31C}, \cite{T32C}, as well as the survey papers  \cite{abikoff},  \cite{ji&papadop}, \cite{acampo&ji&papadop}.

Concerning Teichm\"uller's style, we quote Ahlfors and Gehring, from the preface to Teichm\"uller's \emph{Collected works} \cite{T-collected}: ``Teichm\"uller's style was unorthodox, to say the least. He himself was well aware of the difference between a proof and an intuitive reasoning, but his manner of presentation makes it difficult to follow the frequent shifts from one mode to another." One can also be reminded here of the fact, recalled by Abikoff in \cite{abikoff}, that the style of the journal \textit{Deutsche Mathematik} in which most of Teichm\"uller's articles appeared was to concentrate on general ideas and not on technical details. This was supposed to be a certain idea of ``German mathematics," in the tradition of Riemann and Klein; see \cite{T31C}.  Abikoff in \cite{abikoff} writes: ``The tradition of \emph{Deutsche Mathematik} is one of heuristic argument and contempt for formal proof." He adds that he learned from a conversation with Herbert Busemann that  ``Teichm\"uller manifested those traits early in his career but when pressed could offer a formal proof."

Several results in the paper \cite{T20} are given with only sketchs of proofs, and are stated as ``conjectures."\footnote{The word ``conjecture," in this setting, sometimes means a claim which is not proved immediately, but which is proved later in the paper. For instance, a conjecture is made in \S\,100, and  \S\,101  starts with: ``This extremely insufficiently grounded conjecture shall now be proved."} Regarding this fact, let us quote Teichm\"uller, from the introduction of his later paper \cite{T29}, where he talks about the paper which we are commenting here: ``In 1939, it was a risk to publish a lengthy article entirely built on conjectures. I had studied the topic thoroughly, was convinced of the truth of my conjectures and I did not want to keep back from the public the beautiful connections and perspectives that I had arrived at. Moreover, I wanted to encourage attempts for proofs. I acknowledge the reproaches, that have been made to me from various sides, even today, as justifiable but only in the sense that an {\it unscrupulous\/} imitation of my procedure would certainly lead to a barbarization of our mathematical literature. But I never had any doubts about the correctness of my article, and I am glad now to be able to actually prove its main part."

The setting of Teichm\"uller's paper is that of general surfaces,  orientable or not, with boundary,  or without boundary and with or without distinguished points in the interior or on the boundary.

We have divided the commentary into sections which follow the ordering of  the various sections of the paper \cite{T20} (although the sections of the paper \cite{T20} are not numbered).

Chapter 17 of Volume VI of the present Handbook contains comments by Gr\"otzsch on this paper of Teichm\"uller, edited by R. K\"uhnau \cite{Kuehnau-H}.

\medskip

\noindent{\it Acknowledgements.---}The authors would like to thank Melkana Brakalova-Trevithick 
for reading the manuscript and making corrections.

\section{Introduction}

Teichm\"uller's introduction is unusual in the sense that it contains no explicit statement of any result, and no outline of the paper. At the beginning of this introduction, the author says that he will examine ``the behavior of conformal invariants under quasiconformal mappings" and that this leads to the problem of ``finding the maps that deviate the least from conformality under certain additional conditions."  We recall in this respect that a quasiconformal mapping in the sense of Teichm\"uller is 
a one-to-one mapping (he adds: with certain differentiablity properties) whose dilatation quotient is bounded. The dilatation quotient, at a point where the mapping is continuously differentiable, is defined as the ratio of the  big to the small axis of an infinitesimal ellipse which is the image of an infinitesimal circle at that point. The dilatation quotient of the mapping is the supremum over the surface of these dilatation quotients (whenever they are defined).
Teichm\"uller gives a formula for the dilatation quotient in terms of the partial derivaties of the mapping at a point where it is continuously differentiable. More precisely, for such a mapping from the $(x,y)$-plane to the $(u,v)$-plane, the quasiconformal ratio $D$ is given by
\[D=\vert K\vert+\sqrt{K^2-1}=e^{\arccos \vert K\vert},\]
where
 \[
K=\frac{1}{2}\frac{u_{x}^{2}+u_{y}^{2}+v_{x}^{2}+v_{y}^{2}}{u_{x}v_{y}-v_{x}u_{y}}.\]

Teichm\"uller announces a solution of the problem of finding extremal mappings, that is, mappings that deviate the least from conformality, but he declares that he will not be able to provide rigorous proofs. He then writes that this solution of the problem will use the notion of a quadratic differential\index{quadratic differential!meromorphic}. We recall that a meromorphic quadratic differential\index{quadratic differential!meromorphic}\index{meromorphic quadratic differential} on a Riemann surface is an object which is written in a local coordinate $z$ as $Q(z)dz^2$, where $Q$ is a meromorphic function, and where under a change of coordinates, from a coordinate $z$ to a coordinate $z'$, the corresponding functions $Q(z)$ and $\tilde{Q}(z')$ satisfy:
\[\tilde{Q}(z')=Q(z)\left(\frac{dz}{dz'}\right)^2.\]

 The relation between meromorphic quadratic differentials and the  extremal problem under consideration is indeed one of the major elements in this paper. Teichm\"uller associates a  quadratic differential to an extremal problem between two Riemann surfaces. In particular, Teichm\"uller shows that the trajectory structure of this quadratic differential gives the direction of stretching of the best quasiconformal map between the two surfaces. The metric on the Riemann surface induced by the quadratic differential also plays an important role in Teichm\"uller's solution of the extremal problem.

The rest of the introduction is very broad, without being precise. The author essentially says that he will ``compute the number of conformal invariants" by using the Riemann-Roch theorem, and he mentions the scope of the tools involved:  from computations with infinitely small quantities to uniformization theory, from the length-area method to the domain invariance theorem, from integration of partial differential equations to Galois theory. He warns the reader over and over again that the main result will not be proved, and he justifies this fact by saying that ``proving means reversing the train of thought." This expression will be used once more in \S\,115. At the end of the introduction, it is impossible for the reader to know with precision what is proved in the paper and what is sketched without proof.  As we already mentioned, the details of the proof of one of the main results, the famous \emph{Teichm\"uller Existence Theorem},\index{Teichm\"uller existence theorem}\index{theorem!Teichm\"uller existence} will be provided a few years later, in \cite{T29}. 

\section{Examples of conformal invariants}\label{section2}

The title of this section is explicit enough. The author exhibits complete systems of conformal invariants for some special surfaces. In \S\,2, he considers doubly connected domains in the plane (surfaces homeomorphic to annuli). Such a surface has only one conformal invariant, the modulus of the annulus that is conformally equivalent to it. He gives the group of biholomorphisms of this domain. This group has one continuous parameter. It is generated by rotations and an inversion. In \S\,3, Teichm\"uller considers three examples of surfaces. The first example is a simply connected domain. Such a domain has no conformal invariant, since by the Riemann Mapping Theorem,\index{theorem!Riemann mapping}\index{Riemann mapping theorem} any simply connected domain, provided it is not the entire complex plane, can be mapped conformally onto the unit disk. The second example is a simply connected domain  with one distinguished inner point. As in the previous case, such a domain has no conformal invariant. The third example is a simply connected domain with two  distinguished points in the interior. It has one conformal invariant, the Green function.\index{Green function} In \S\,4, Teichm\"uller first considers the case of a triangle, i.e. a simply connected domain with three distinguished points on the boundary. He recalls that such a domain has no conformal invariant. A quadrilateral, that is, a simply connected domain with four (ordered) distinguished points on the boundary, has one conformal invariant, which Teichm\"uller calls a \emph{characteristic} conformal invariant. He defines it as the cross ratio of the four ordered points. It is the number $\lambda$ obtained by mapping the domain conformally onto the upper half-space while mapping the four points in the given order to the points $0, 1,\lambda,\infty$.  In \S\,5, he considers the case of the torus.  This subsection is interesting because the notion of \textit{Teichm\"uller space} (in the usual sense) appears for the first time. Indeed, Teichm\"uller recalls that a torus is determined by a \emph{period}, an element of the upper half-plane $\mathbb{H}$, and that two tori are biholomorphic if and only if their periods differ by an element of $\mathrm{SL}(2,\mathbb{Z})$. He recalls that the conformal invariant in this case is the modular function of the given complex torus. He then says that it is better to consider the period instead of the modular function. A marking by a pair of simple closed curves is also considered. This is called a \emph{canonical cut} in \S\,50. Thus, he considers the Teichm\"uller space of the torus. He proves that this space is the hyperbolic plane $\mathbb{H}$.\index{Teichm\"uller space!torus}\index{torus!Teichm\"uller space}\index{torus!Teichm\"uller space!metric} We also note that it is the first time where the term ``deformation" appears. In \S\,6, Teichm\"uller considers a doubly connected domain with two distinguished points on the boundary. He shows that in this case there are exactly two conformal invariants, distinguishing between the cases where the two distinguished points  belong to the same boundary component or not.

In conclusion, this section is a kind of introduction to the notion of conformal invariant. The Teichm\"uller space of the torus will be more thoroughly examined in Sections \ref{sectiontore1}  and 29. Conformal invariants will be used later on in establishing the ``dimension formula," Equation (\ref{dimension}) below.

\section{Non-orientable regions}

In \S\,7 and \S\,8, Teichm\"uller defines, through two examples, the notion of conformal invariant for ``non-orientable regions." Such an object is invariant by a conformal mapping. In this respect it is important to note that according to Teichm\"uller's definition, a conformal mapping (even in the orientable case) preserves only the angles, but not necessarily the orientation. Such a mapping is also defined in the non-orientabe case. In \S\,7, he considers the case of the M\"obius band. The conformal invariant of such a surface is the outer radius of its 2-sheeted unbranched covering, which is an annulus. In \S\,8, he considers the projective plane  $\mathbb{R}\textrm{P}^2$ (which he calls the \textit{elliptic plane}) with two distinguished points. Using the sphere in $\mathbb{R}^3$ as a 2-sheeted covering of this surface, he shows that there exists in this case only one conformal invariant, the shortest length (with respect to the spherical metric) between the corresponding distinguished points on the sphere.

Thus, in these two examples of non-orientable surfaces, the problem of finding the conformal invariants is lifted to the oriented double cover. This idea will be used to solve the main problem of extremality (Problem \ref{probexistence&unicite} below); see \S\,21 and 25.

\section{Finite Riemann surfaces}

In \S\,9 and \S\,10, Teichm\"uller recalls the notion of a Riemann surface. He uses the expression ``Riemannische Mannigfaltkeit" for such a surface.\index{Riemann surface}  He gives the classical definition of a surface, as a two-dimensional manifold which admits a triangulation. In accordance with the convention already mentioned, the transition maps are not necessarily holomorphic but they may be anti-holomorphic. Teichm\"uller considers surfaces of topological finite type (i.e. whose fundamental group is finitely generated), and he calls them ``finite Riemann surfaces."\index{Riemann surface!finite} They may be bordered or not, planar or not, orientable or not. Surfaces of infinite topological type (those whose  fundamental group is infinitely generated) are not considered.  Note that according to Teichm\"uller's convention, a ``finite Riemann surface"  has no distinguished points. Surfaces with distinguished points will be considered in the next section, and they are given a special name.

\section{The notion of a principal region}\label{section5}

Teichm\"uller continues giving definitions, while considering examples. He introduces the notion of  a\index{principal region} \textit{principal region}.\footnote{In \cite{T23}, Teichm\"uller, explains his choice of the name ``principal region," while putting it in a much wider content. He writes: ``I have already encountered in various fields this concept: An object $\mathfrak{A}$ (here, a principal region) which is only determined if an object $\mathfrak{A}_1$ (here, a manifold) is first given, and then an object $\mathfrak{A}_2$ (here, the distinguished points) is defined using $\mathfrak{A}_1$. Similarly: in algebra, $\mathfrak{A}_1$, a cyclic field extension of a fixed field,  $\mathfrak{A}_2$: generators of the Galois group, $\mathfrak{A}$: cyclic normal field." He refers to this example in \cite{T4}. In the same paper, \cite{T23}, Teichm\"uller introduces the notion of map between principal regions, which is adapted to this general setting.}\index{principal region} This is a finite Riemann surface in which $h\geq 0$ inner points and $k\geq 0$ boundary points are distinguished.  He defines the notion of a mapping between principal regions that take distinguished points to distinguished points.
He mentions that surfaces with punctures are not considered because they can be treated as surfaces with distinguished inner points.

Teichm\"uller announces in \S\,12 a forthcoming work in which distinguished points are not fixed but are moved infinitely close to each other. This is an idea of degeneration of Riemann surfaces. In \S\,13, he denotes by $\mathfrak{R}^\sigma$ the space we call today the\index{Riemann moduli space} \textit{Riemann moduli space}.\footnote{Note the use of the supersript $\sigma$ in $\mathfrak{R}^\sigma$, instead of indexing by the genus, number of boundary components, orientability, etc. as we usually do today. Teichm\"uller's notation highlights the dimension, which is an important element in his work, but it does not make a distinction between the various moduli (and Teichm\"uller) spaces of equal dimensions.}  The exponent $\sigma$ refers here to the number of conformal invariants. Teichm\"uller writes about this space:
\begin{quote}
We shall admit the following statement without proof: when one identifies conformally equivalent principal regions, principal regions of fixed topological type form a topological manifold which is on a small scale homeomorphic to the $\sigma$-dimensional Euclidean space and hence can be described as the space $\mathfrak{R}^\sigma$.

The conformal invariants of the principal region are then precisely the functions on $\mathfrak{R}^\sigma$. Hence there are,  on a small scale, precisely $\sigma$ independent conformal invariants.
\end{quote}
It is known that the moduli space has the structure of an orbifold,\footnote{This is Thurston's terminology.} that is, a space with a special kind of local singularities (quotiens of actions of finite groups). The singularities are due to the existence of finite order elements in the mapping class group acting on Teichm\"uller space. In a subsequent paper \cite{T32}, Teichm\"uller points out explicitly that there are singularities in $\mathfrak{R}^\sigma$.\footnote{Cf. the English translation of \cite{T32}, in Volume IV of this Handbook, p. 788: ``It turns out that $\mathfrak{R}$ contains certain singular manifolds. But we will construct a
covering space $\underline{\mathfrak{R}}$ without singularities."} After Section 15 of the paper under consideration, Teichm\"uller will consider a covering space of $\mathfrak{R}^\sigma$,  denoted by $R^\sigma$, which corresponds to the Teichm\"uller space in the case of oriented surfaces. This is a covering of $\mathfrak{R}^\sigma$ with no singularities.\footnote{We recall by the way that Riemann's  moduli space has a \emph{finite} covering which is non-singular, or, equivalently, that the mapping class group has a torsion-free subgroup of finite index.}

\section{Characteristic numbers}\label{dimension}

Teichm\"uller starts this section by giving the \textit{dimension formula} for moduli space. The proof of this formula is postponed to \S\,112. The formula is the following:
\begin{equation}\label{eq1}
\sigma - \rho = 6g-6+3\gamma+3n+2h+k.
\end{equation}
The left hand side denotes conformal objects, and the right hand side topological objects. Precisely, the letters denote the following:
\begin{itemize}
\item $g$ is the number of handles,
\item $\gamma$ is the number of crosscaps,
\item $n$ is the number of boundary curves,
\item $h$ is the number of distinguished inner points,
\item $k$ is the number of distinguished boundary points,
\item $\rho$ is the number of parameters of the (countinuous) group of conformal mappings of the principal region onto itself,
\item $\sigma$ is the dimension of the space of all classes of conformally equivalent  principal regions which are topologically equivalent to the given principal region. This is Riemann's ``number of moduli."
\end{itemize}
Teichm\"uller recalls the notion of \emph{algebraic genus}\index{genus!algebraic}, and  introduces the notion of  \emph{reduced dimension}.\index{reduced dimension}\index{dimension!reduced} These notions take into account the non-orientability of the surface. The algebraic genus of a surface was defined by Klein in his study of automorphisms of orientable and non-orientable Riemann surfaces, with or without boundary. For a non-orientable surface, the algebraic genus is the genus of its orientable double cover. The reduced dimension of a non-orientable surface or a surface with boundary is also defined in terms of a double cover. The result that Teichm\"uller states is that for bordered or non-orientable surface, the reduced dimension $\tau$ is equal to $\sigma$, and for  a closed orientable surface, the reduced dimension $\tau$ is equal to $\sigma/2$.

The notions of  algebraic genus  and reduced dimension play an important role in the proof of Equation (\ref{formuledimension}) in the case of non-orientable surfaces or bordered surfaces. To prove Formula (\ref{eq1}), Teichm\"uller uses the Riemann-Roch theorem and shows that $\sigma$ (resp. $\rho$) corresponds to the dimension of the subspace of the space of meromorphic quadratic differentials\index{inverse differential} (resp. inverse differentials).\footnote{``Inverse differentials" are also called ``reciprocal differentials." The German word is ``reziproken Differenziale." We chose the first translation, which is the one used by Ahlfors.} See Sections 13 and 21 for the corresponding definitions. An inverse differential on a Riemann surfaces is a  section of the tangent bundle, that is, a vecor field. 

We note that Formula (\ref{eq1}) is more elegant than the one that is usually given today (in which one makes a distinction between the cases of the torus and the other cases) since it applies to any genus $\geq 0$ and to any number of boundary components.  More importantly, the formula, as it is stated here, bears some strong resemblance to the Riemann-Roch formula, and in fact, it may be considered as a ``global Riemann-Roch formula"\index{theorem!Riemann-Roch (global)} where, on the left hand side, $\sigma$ is the dimension of a moduli space which is a global version of the vector space of quadratic differentials, and $\rho$ is the dimension of a continuous group of conformal maps of the surface, which is a global version of the vector space of tangent vectors. The right hand side involves topological data.

\section{The statement of the problem}\label{sectionex}

This section is important for the rest of the paper. Teichm\"uller uses quasiconformal mappings to define what we call today the  \textit{Teichm\"uller distance}.\index{Teichm\"uller metric}\index{metric!Teichm\"uller} He also states that this distance is Finsler. At this point, he does not introduce Teichm\"uller space, but  he rather continues to work on Riemann's moduli space. Let us give some more details.

In \S\,15, an analytic formula is given for the dilatation quotient\index{dilatation quotient} at a point for maps between two domains of the complex plane. The result is said to be equal to the ratio of the great axis to the small axis of an infinitesimal ellipse which is the image of an infinitesimal circle at the given point. The logarithm of the dilatation quotient is a measure of the deviation of the map from being conformal.

In \S\,16, considering conformal invariants as functions from $\mathfrak{R}^\sigma$ to the reals, Teichm\"uller formulates the following problem:

\begin{problem}\label{prob1}
Consider a fixed principal region. Let $J$ be a conformal invariant considered as a function on $\mathfrak{R}^\sigma$, and let  $C>1$ be given. What values does $J$ assume for those principal regions onto which the given principal region can be quasiconformally mapped so that the dilatation quotient is everywhere $\leq C$?
\end{problem}

Thus, instead of searching for conformal invariants, Teichm\"uller searchs now for the behaviour of invariants under quasiconformal mapppings of bounded dilatation quotient, and for  precise estimates on this behavior in terms of this dilatation quotient. He then states  the following equivalent form of Problem \ref{prob1}:

\begin{problem}\label{prob2}
Let $P$ be an element of $\mathfrak{R}^\sigma$ and  $C>1$ a given number. We map in all possible ways the principal region represented by $P$ quasiconformally onto another principal region of the same type so that the upper bound of the dilatation quotients is $\leq C$, respectively $<C$, and we represent the obtained principal region by a point $Q$ of $\mathfrak{R}^\sigma$. Give a description of the set $\overline{\mathfrak{U}_{C}\left( P \right)}$, respectively $\mathfrak{U}_{C}\left( P \right)$, of the  points $Q$ so obtained.
\end{problem}

The definition of the set $\mathfrak{U}_{C}\left( P \right)$ is given after the introduction of a distance (the Teichm\"uller distance), in \S\,18. The precise definition is the following:

\begin{definition}
Let two principal regions be given, represented by two points $P$ and $Q$ of $\mathfrak{R}^\sigma$. For a quasiconformal mapping from this principal region onto another, let $C$ be the upper bound on the surface of the dilatation quotients. We set the logarithm of the lower bound of all these $C$ to be  the distance $\left[ PQ \right]$ between the two points or between the principal regions.
\end{definition}

Using this definition, the set $\mathfrak{U}_{C}$ is:
\[
\mathfrak{U}_C \left( P \right)= \left\lbrace Q\in\mathfrak{R}^\sigma \; \mid \; \left[ PQ \right] < C \right\rbrace.
\]
Teichm\"uller shows that $[\cdot \cdot]$ is a distance.\footnote{We note that the separation property for this metric, that is, the fact that if the quasiconformal constant is equal to 1 then the map us conformal, in the more general setting where the quasiconformal dilatation is a distribution, is a corollary of the so-called Weyl Lemma.} The collection of sets
 $$\left\lbrace \mathfrak{U}_{C}\left( P \right) \, \mid \; C\geq 0, \; P\in\mathfrak{R}^\sigma \right\rbrace$$
 forms a  system of neighborhoods which equips $\mathfrak{R}^\sigma$ with a topology.

In \S\,19, Teichm\"uller claims that he will show that $\left( \mathfrak{R}^\sigma , [..] \right)$ is Finsler. In fact, he will show that the metric space $\left( R^\sigma , \left[ \cdot  \cdot \right] \right)$ is Finsler.  (Recall that $R^\sigma$ is Teichm\"uller space.) He proposes the following problem, which is analogous to Problem \ref{prob2}, but which concerns mappings instead of spaces.
\begin{problem}\label{prob3} Given two planar principal regions represented through the points $P$ and $Q$ of $\mathfrak{R}^\sigma$, determine all the extremal quasiconformal mappings whose dilatation quotient is everywhere $\leq e^{[PQ]}$.\end{problem}

An extremal mapping will be a map whose quasiconformal dilatation attains the infimum. Problem \ref{prob3} suggests the study of the following two questions:

\begin{itemize}
\item the existence of extremal mappings;
\item the description of a geodesic segment between $P$ and $Q$.
\end{itemize}

The aim of the next four sections is to solve Problem \ref{prob3} in particular cases.
\section{Examples of bordered orientable principal regions}\label{section8}

In this section, Teichm\"uller gives explicitly the Teichm\"uller distance and the extremal mapping between some special surfaces, namely, those studied in Section 2. For that, he uses the conformal invariants that he already determined in Section \ref{section2}.

In \S\,20, he deals with quadrilaterals. He gives a proof of Gr\"otzsch's\index{theorem!Gr\"otzsch} theorem\footnote{Camillo Herbert Gr\"otzsch\index{Gr\"otzsch, Herbert} (1902-1993) is the founder of the theory of quasiconformal mappings. (It is admitted that Ahlfors and Lavrentiev are also founders of this theory, but the work of Gr\"otzsch on this subject precedes their works.) Gr\"otzsch did his PhD under Paul Koebe in Leipzig in 1929. He is also well known for his work on graph theory. Reiner K\"uhnau, who did his PhD under Gr\"otzsch, wrote two long reviews on his work, cf. \cite{Kuhnau1} and \cite{Kuhnau2}. Teichm\"uller has always acknowledged the importance of Gr\"otzsch's work in his research involving quasiconformal mappings. In \cite{T23}, where Gr\"otzsch's work is mentioned several times, Teichm\"uller writes: ``Much of what is discussed here is already contained in the works of Gr\"otzsch, but mostly hidden or specialized in typical cases and in a different terminology." We refer the reader to Chapter 17 of Volume VI of this Handbook, by R. K\"uhnau \cite{Kuehnau-H}.} for  quasiconformal mappings between rectangles, cf. \cite{rd4}, which solves Problem \ref{prob3} in that special case. The precise statement of this theorem, stated in Teichm\"uller's paper, is the following:

\begin{theorem}\label{prop1}
If we map the rectangle $0< \mathrm{Re}\left( \zeta \right)<a$, $0<\mathrm{Im}\left( \zeta \right)< b$ onto the rectangle $0< \mathrm{Re}\left( \zeta^\prime \right)<a^\prime$, $0<\mathrm{Im}\left( \zeta^\prime \right)< b^\prime$, such that the vertices $\zeta =0, a ,a+\mathrm{i} b, \mathrm{i} b$ are sent to $\zeta^\prime =0, a^\prime ,a^\prime+\mathrm{i} b^\prime, \mathrm{i} b^\prime$ and the dilatation quotient is everywhere $\leq C$, then $\frac{b}{a}\leq C \frac{b^\prime}{a^\prime}$. Equality holds only for an affine mapping.
\end{theorem}

  The method of proof is usually called
the\index{length-area method} \emph{length-area method}.\footnote{In his \emph{Collected Works} \cite{Ahlfors:Collected} Vol. 1 (p. 1),  talking about his first two published papers on asymptotic values of entire functions of finite order, Ahlfors comments on the origin of this method. He writes: ``The salient feature of the proof is the use of what is now called the length-area method. The early history of this method is obscure, but I knew it from and was inspired by its application in the well-known textbook of Hurwitz-Courant to the boundary correspondence in conformal mapping. None of us [Ahlfors is talking about Nevanlinna and himself] was aware that only months earlier H. Gr\"otzsch had published two important papers on extremal problems in conformal mapping in which the same method is used in a more sophisticated manner. My only priority, if I can claim one, is to have used the method on a problem that is not originally stated in terms of conformal mapping. The method that Gr\"otzsch and I used is a precursor of the method of extremal length [...]  In my thesis \cite{Ahlfors:Unter} the lemma on conformal mapping has become the main theorem in the form of a strong and explicit inequality or distorsion theorem for the conformal mapping from a general strip domain to a parallel strip, together with a weaker inequality in the opposite direction. [...] A more precise form of the first inequality was later given by O. Teichm\"uller." (Ahlfors refers here to Teichm\"uller's paper \cite{T200}.)}   

The result implies that for two quadrilaterals $P$ and $Q$ of moduli $m_P$ and $m_Q$ respectively, we have
$$
[PQ]=\vert \log\left( m_P \right) - \log\left( m_Q \right) \vert
$$ñ
and also
$$
\mathfrak{U}_{C}\left( P \right)= \left\lbrace Q \; \mid \; \frac{1}{C}m_P <m_Q <C m_P  \right\rbrace.
$$
This also gives, using Teichm\"uller's notation, that  $\mathfrak{R}^1 = \mathbb{R}^{*}_{+}$.

In \S\,21, Teichm\"uller solves the same problem for a doubly connected domain  (called a \textit{ring domain}\index{ring domain} in the translation). As we already noted, the conformal type of such a domain is determined by one parameter called the modulus $M$. Teichm\"uller proves the following.
\begin{proposition}
For a quasiconformal mapping of a ring domain with modulus  $M$ onto a ring domain with modulus  $M^\prime$, we have $M^\prime \leq C M$.
\end{proposition}

 Here, the value $C$ is the upper bound of the dilatation quotient. Teichm\"uller  also gives the extremal mapping (therefore the Teichm\"uller distance) and the description of the moduli space of annuli, $\mathfrak{R}^1$. Like the moduli space of quadrilaterals, the moduli space of annuli can be identified with $\mathbb{R}^{*}_{+}$. This is not a surprise. Indeed, in \S\,20, Teichm\"uller shows that the moduli problem for annuli is equivalent to the one for quadrilaterals.

In \S\,23, Teichm\"uller determines the moduli space of simply connected domains with two distinguished points in their interior. He also computes the Teichm\"uller distance between two such domains. One interesting feature of the solution is the use of ramified covering theory, a device that allows him to reduce the question to a question already solved, namely, that of moduli of doubly connected domains. We shall see that this use of covering theory is much more general. The Gr\"otzsch ring function\index{Gr\"otzsch function} $\mu(r)$ appears in this section. This function (the modulus of the ring domain obtained by cutting the unit disk along the segment $[0,r]$ ($0<r<1$), is considered in Gr\"otzsch's paper \cite{Gr} and Teichm\"uller's paper \cite{T200}). See also Teichm\"uller's paper \cite{T31} and the commentary \cite{T31C} in the same volume for a study of such a domain. The Gr\"otzsch function plays  an  important role in the development of the theory of quasiconformal mappings.

In \S\,24, Teichm\"uller considers simply connected domains with two distinguished points on their boundary and one in their interior. In mapping conformally such a domain onto a disk, the interior point can be sent  to $0$ and one of the two boundary points to $1$. The image of the second boundary point is a conformal invariant of such a domain. Like in the previous case, Teichm\"uller uses the covering map $z\mapsto z^2$ to reduce the problem to a question already solved, namely,  that of moduli of quadrilaterals. He determines the extremal mapping (and therefore the Teichm\"uller distance)  and the corresponding moduli space. The idea of using the covering  $z\mapsto z^2$ to solve an extremal problem is also used in Teichm\"uller's subsequent paper \cite{T31}; see also the commentary \cite{T31C}.

To conclude this section, let us make a few remarks on these examples. In all these cases, the moduli space is identified with $\mathbb{R}^{*}_{+}$. This agrees with the dimension count given in Section \ref{dimension}, since the automorphism groups of the corresponding regions are trivial. We also remark that the Teichm\"uller distance is shown to be equal to the absolute value of the logarithm of ratios of numbers depending on corresponding conformal invariants. This fact is related to the \textit{Kerckhoff formula}\index{Kerckhoff formula} which gives the Teichm\"uller distance as a supremum of ratios of extremal lengths(see \cite{kerckhoff}). Finally, in the last two  subsections (\S\,23 and \S\,24), there appear functions which play the role of quadratic differentials. These functions are used to solve directly the problem, i.e without using covering theory.

\section{The torus}\label{sectiontore1}

In this section, Teichm\"uller solves completely Problem \ref{prob3} in the case where the surface is a torus. In fact,  he solves this problem by considering Teichm\"uller space and not the moduli space. He also shows that all the previous examples can be deduced from the case of the torus.

In \S\,25, he considers the two primitive periods $w_1$ and $w_2$ associated with a torus, and their quotient ratio $\omega=\frac{w_1}{w_2}$. He calls this quotient (which can be assumed to have positive imaginary part) the ``characteristic conformal invariant" of the torus. By passing to the universal covering space, he shows:

\begin{proposition}
Let $T_1$ and $T_2$ be two tori of period ratios $\omega_1$ and $\omega_2$ respectively. Let $f: T_1 \rightarrow T_2$ be a quasiconformal mapping whose quasiconformal dilatation is bounded by $C$. Then
$$
\frac{1}{C}\mathrm{Im}\left( \omega_1 \right) \leq \mathrm{Im}\left( \omega_2 \right) \leq C \mathrm{Im}\left( \omega_1 \right),
$$
with equality if and only if its lift $\tilde{f}: \mathbb{C} \rightarrow \mathbb{C}$ is an affine map. 
\end{proposition}

In \S\,26, he shows:
\begin{proposition}
\textit{The lower limit for the logarithms of all upper bounds for dilatation quotients of the mapping from a torus with period ratio $\omega$ onto a torus with period ratio $\omega^\prime$ is equal to the non-Euclidean distance between $\omega$ and $\omega^\prime$ in the upper half-plane.}
\end{proposition}
Thus, the Teichm\"uller metric of the torus\index{Teichm\"uller space!torus}\index{torus!Teichm\"uller space}\index{torus!Teichm\"uller space!metric} is the hyperbolic metric on the upper half-plane. This was already stated in Section \ref{section2}. The present section contains new material.

In \S\,27, he shows that Poblem \ref{prob3} for the sphere with four distinguished points can be reduced to the case of the torus. This is done by taking the two-sheeted branched covering with branching at the four distinguished points. He shows that the moduli space of such a space with respect to the Teichm\"uller metric is identified with unit disc equipped with its hyperbolic metric. He concludes by saying that the moduli space is a \textit{complete differential geometric space}.

In \S\,28, he proves what we call now the  \textit{Teichm\"uller existence and uniqueness theorems}, in the case of the torus.\index{Teichm\"uller existence theorem!torus}\index{Teichm\"uller uniqueness theorem!torus}

In \S\,29, he shows that all the examples from Section \ref{section8} can be deduced from the torus case. In fact, it suffices to consider the annulus, since all other examples are deduced from that case. For this, he uses a doubling process  along two boundary components, which gives the torus. This doubling process will be thoroughly used in Section  \ref{sectionduplication}.

\section{Examples of non-orientable principal regions}\label{sectionexnonorientable}

In this section Teichm\"uller considers three different non-orientable surfaces: the M\"obius band, the projective plane with two distinguished points, and the Klein bottle. He solves, for these surfaces, the two problems that he solved in the preceding two sections for the special orientable surfaces, that is, the problem of finding the conformal invariants and the problem stated in Section \ref{sectionex}.

In \S\,30, he considers the M\"obius band. By using its two-sheeted covering (an annulus), he gets the Teichm\"uller distance.

In \S\,31, he considers the projective plane with two distinguished points. He shows that this case can be reduced to the case of the sphere with four distinguished points.

In \S\,32, he considers the Klein bottle and he lifts the problem to a problem concerning the case of the torus. The torus is the two-sheeted orientation covering of the Klein bottle.

These three examples are special cases of the fact that one deals with non-orientable surfaces by passing to the orientation double cover. This will be thoroughly used in Section \ref{sectionduplication}.

\section{A wrong track}
In \S\,33, Teichm\"uller explains why the method used in the preceding sections does not work in the case of an annulus with two distinguished points belonging to different boundary components.  He will solve the problem in this case in \S\,131 (Section 26).

In \S\,34, he gives some hints to how to generalize the above examples. For example, he says that  for  closed oriented surfaces, we have to consider their universal covering which is in this case the non-euclidean plane. This is where hyperbolic geometry  enters into this theory.

The result of \S\,35 can be considered (using today's terminology) as a comparison between the Teichm\"uller distance and the so-called \textit{Thurston metric}. The author establishes an inequality on the effect of a quasiconformal mapping conjugating two groups of hyperbolic transformations of the non-Euclidean plane. He obtains an inequality involving the dilatation quotient of the quasiconformal mapping and the dilatations of corresponding hyperbolic transformations, or, equivalently, the lengths of corresponding closed geodesics in the quotient surfaces. In this way, Teichm\"uller proves what we call today the \textit{Wolpert inequality}.\index{Wolpert inequality} The proof is the same as in \cite{sorvali} and \cite{wolpert} where the authors independently redicovered this formula. A question on the topologies defined by these two metrics is addressed, and Teichm\"uller says: ``I am not able to answer this question."

\section{Extremal quasiconformal mappings}

In this short section, some important notions are introduced in a heuristic way, namely, the two foliations associated with   a quadratic differential (the horizontal and the vertical foliation) and the notion of Teichm\"uller geodesic.

In \S\,36, Teichm\"uller argues that the quasiconformal dilatation of an extremal mapping is constant. He says however that ``this should not be any proof whatsoever, rather only a heuristic consideration."

In \S\,37, he asks explicitly: ``which direction fields are associated with extremal quasiconformal mappings?"

In \S\,38, he defines the notion of \textit{locally extremal quasiconformal mapping}\index{locally extremal quasiconformal mapping}\index{quasiconformal mapping!locally extremal} as follows:
\begin{definition}\label{def1}
A quasiconformal mapping is called locally extremal if its dilatation quotient is a constant $=K$ and if there is, for each point $\mathfrak{p}$ in the open set where the mapping is defined, a simply connected neighborhood $\mathfrak{U}$ with the following property: the mapping is defined in $\mathfrak{U}$ and on the boundary of $\mathfrak{U}$ and each mapping which agrees in the complement of $\mathfrak{U}$ with the given one and whose dilatation quotient is $\leq K$, has constant dilatation quotient $K$.
\end{definition}

In \S\,39, Teichm\"uller introduces the notion of geodesic segment with respect to the Teichm\"uller distance. This is obtained by using the same pair of transverse foliations on the surface, and varying the dilatation constant.

\section{The Riemann-Roch theorem}

This section makes a link, through the space of holomorphic quadratic differentials, between the Riemann-Roch theorem\index{Riemann-Roch theorem}\index{theorem!Riemann-Roch} and the dimension formula for the Riemann moduli space, for closed orientable surfaces.\footnote{In the paper \cite{T23}, a generalized version of the Riemann-Roch theorem is obtained, which includes  special data at the distinguished points.}  Teichm\"uller proves that the complex dimension of the latter is equal to $3g-3$. In \S\,40, he rewrites the dimension formula only for closed surfaces. Using his  notation, this formula is:
$$
\sigma -\rho =6\left( g-1\right).
$$
In \S\,41, he recalls the notion of \textit{divisor},\index{divisor} \textit{prime divisor}\index{prime divisor}\index{divisor!prime} (associated with a point) and \textit{principal divisor}\index{principal divisor}\index{divisor!principal} (associated with a meromorphic function) on a closed surface. He uses the multiplicative notation for the combination of divisors	. Prime divisors generate a free abelian group. He  introduces the notion of \textit{divisor class},\index{divisor class} an element of the \textit{divisor class group}, the group of divisors quotiented by the normal subgroup of principal divisors.

In \S\,42, he introduces the notion of \textit{n-dimensional differential $d\zeta^n$}, describing the behaviour with respect to coordinate change which makes this notion invariant. In the holomorphic local coordinate $t$,  an $n$-dimensional differential $d\zeta^n$ has the form $d\zeta^n=g(t)dt^n$ where $g(t)$ is a (local) holomorphic function of $t$. The invariance property says that if $t'$ is another holomorphic local coordinate, and if the differential has the form $g'(t)dt'^n$ in this coordinate, then at the overlap of the coordinates we have
\[g(t)=g'(t')\left(\frac{dt'}{dt}\right)^n
.\]
For $n=1$, one gets the usual notion of (meromorphic) differential form, and for $n=2$ the notion of quadratic differential.\index{quadratic differential} The behavior at the distinguished points will be specified.

To each $n$-dimensional differential is associated a divisor. The space of holomorphic $n$-differentials is denoted by $\mathfrak{M}^n$.

In \S\,43, Teichm\"uller defines the notions of \textit{degree}\index{degree!divisor}\index{divisor!degree} of a divisor and of \textit{integral divisor}. This is needed in the statement of the Riemann-Roch theorem, which he applies  to the case of $n$-dimensional differentials.  
In \S\,44, he considers the $2$-dimensional differentials, i.e quadratic differentials and $(-1)$-dimensional differentials. He calls the latter \textit{inverse differentials}.\index{inverse differential} 
He gives the (complex) dimension of the spaces of such differentials, which is equal to $\sigma$ and $\rho$ respectively.
Thus, using his notation, we have, for closed orientable surfaces without distinguished points,
 \[2\ \mathrm{dim}\  (\mathfrak{M}^2)=\sigma\]
and
\[2\ \mathrm{dim} \left(\frac{1}{\mathfrak{M}}\right)=\rho.\]

Teichm\"uller says that he first noticed these dimension formulae for very particular cases and then he conjectured the rest to get the general dimension formula: 
``As soon as I knew the values of $\sigma$ given in \S\,40, I noticed that they were also equal to $2\textrm{dim} (\mathfrak{M}^2)$. Then, I conjectured the relationship; 
this is what the following discussion is going to be about."

Thus, it is from meditating upon the Riemann-Roch theorem that Teichm\"uller obtained the idea of the dimension formula for Teichm\"uller space. The formula will be proved in Section 21.

\section{A conjecture}\label{sectionconj}

 In this section, Teichm\"uller states the relation between the space of quadratic differentials and the value $\sigma$, and he proves that $\rho$ is equal to twice the complex dimension of the space of inverse differentials. From the latter, he announces what we call today the \textit{Teichm\"uller theorem}\index{theorem!Teichm\"uller}\index{Teichm\"uller theorem} (see the statement in our comments on \S\,46 below).

In \S\,45, the notion of \emph{infinitesimally conformal self-mapping}\index{infinitesimally conformal mapping} of a Riemann surface $\mathcal{F}$ is introduced. This is a mapping which ``moves every point only infinitesimally." The author starts by noting that under an infinitesimally conformal self-mapping,  a point with local holomorphic coordinates $z$ moves to a point with coordinates $z+\epsilon g(z)$, where $\epsilon$ is a constant which is an ``infinitely small quantity," and $g(z)$ is a \emph{field of inverse differentials}, that is, a vector field.

After explaining why $\rho$ is equal to ``the maximal number of (real) linearly independent \textit{infinitesimally conformal mappings}," he shows the following:
\begin{proposition}
Let  $\mathfrak{F}$ be a closed surface.
 The infinitesimally conformal mappings of $\mathfrak{F}$ onto itself correspond bijectively to everywhere finite inverse differentials.
\end{proposition}

Teichm\"uller's next purpose is to investigate the question of which vector fields correspond to extremal infinitesimal conformal mappings. 

In \S\,46, he conjectures the relation between extremal quasiconformal mappings and quadratic differentials. He recalls that an extremal quasiconformal mapping provides a constant $K$ and a direction field and he states the following, which he calls ``the conjecture": 
\begin{theorem}\label{conj1}
 Let $d\zeta^2$ be a nonzero everywhere finite quadratic differential on $\mathfrak{F}$. Let us assign to every point of $\mathfrak{F}$ the direction where $d\zeta^2$ is positive. Then, all extremal quasiconformal mappings are described through the direction fields thus obtained and through arbitrary constant dilatation quotients $K\geq1$.
 \end{theorem}
This conjecture, which is a form of Teichm\"uller's theorem,\index{theorem!Teichm\"uller} will be proved later. It involves the notion of \textit{vertical foliation}\index{vertical foliation} associated with a quadratic differential. For a neighborhood of a non-zero point of $d\zeta^2$, the author introduces \textit{normal coordinates} and the local expression of the \textit{Teichm\"uller mapping}.\index{Teichm\"uller extremal mapping} He proves that  this mapping is locally extremal in the sense of Definition \ref{def1} and admits (temporarily)  that it is also an extremal quasiconformal mapping. He announces that he will prove this in Section 27. He will be concerned with the proof of Theorem \ref{conj1} throughout the rest of the paper. Finally, he shows, based on this therem, that $\sigma=2\dim\left(\mathfrak{M}^2 \right)$.

In \S\,47, Teichm\"uller shows how a quadratic differential determines a field of directions. He considers the expression of a quadratic differential in the neighborhood of a zero and shows that around such a point the mapping, previously introduced, (the so-called ``Teichm\"uller mapping")\index{Teichm\"uller extremal mapping} is, again, locally extremal.

In \S\,48, he considers the case where the quadratic differential has poles. Poles are necessarily of order $1$. To get a locally extremal mapping, he notes that one has to consider ``the pole as a distinguished point of the principal region."

In conclusion, given a quadratic differential, Teichm\"uller introduces a mapping (the ``Teichm\"uller mapping") whose  direction field is exactly the vertical foliation of a quadratic differential, and he shows that this map is a locally extremal mapping.

The preceding sections can be considered as a kind of introduction to the results that will follow. So far, several objects have been introduced from a heuristic point of view. The main result of this paper, namely, the so-called Teichm\"uller theorem,\index{theorem!Teichm\"uller}\index{Teichm\"uller theorem} will be treated next.

\section{Topological determination of principal regions}

This section is devoted to the introduction of the \textit{Teichm\"uller space},\index{Teichm\"uller space} which the author calls ``space of topologically determined principal regions." Teichm\"uller states the main problem from this new point of view. We shall quote this problem below.

In \S\,49, he introduces the notion of a \textit{marked surface}. The term \textit{deformation} is defined precisely and corresponds to the notion already used from a heuristic point of view in Section \ref{sectiontore1}. Teichm\"uller denotes by $R^\sigma$ the space of marked surfaces and observes that there is a canonical map from $R^\sigma$ to $\mathfrak{R}^\sigma$. In fact, he explains that
$$
\mathfrak{R}^\sigma = R^\sigma /\textrm{MCG}\left( \mathfrak{h}_0\right)
$$
where $\mathfrak{h}_0$ is a fixed principal region and $\textrm{MCG}\left( \mathfrak{h}_0\right)$ denotes its mapping class group. He will give more details about this group in Section 28.

In \S\,50, he proves, using the new definition of $R^2$, that the Teichm\"uller space of the torus coincides with the space determined in Section \ref{sectiontore1}. In other words, it is the upper half-plane. He writes:
 ``It is not $\mathfrak{R}^2$ ($\sigma=2$!), but rather the covering space $R^2$ that appears  to be  mapped bijectively onto the upper $\omega$-half-space."
He notes that the Teichm\"uller space (or the space of topologically determined principal regions) was already studied, in some special cases, in Sections \ref{sectionex}, \ref{section8} and \ref{sectionexnonorientable}.

In \S\,51, he explains that the problems formulated in Section \ref{sectionex} can be considered as problems concerning the space $R^\sigma$.

In \S\,52, he reformulates Problem \ref{prob3} in the setting of topologically determined principal regions. The problem is now stated as follows:
\begin{problem}\label{probexistence&unicite}
Let principal region of a certain type be fixed.
One shall give a set of quasiconformal mappings, each with constant dilatation quotient, from some principal region of this type onto another.
One shall prove that every mapping of this set is extremal quasiconformal, in the sense that for any mapping between the same principal regions obtained by deforming it, the maximum of its dilatation quotient is bigger, or equal only if the new comparison mapping is obtained by composing the former mapping which comes from our set with a conformal mapping.
Finally, one shall show that any topological mapping from one principal region of the type considered onto another can be deformed to a mapping from this set.
\end{problem}

Thus, the problem now is to associate with any fixed Riemann surface a set of  quasiconformal mappings from that surface onto the other Riemann surfaces of the same type, such that:
\begin{enumerate}
\item each of these maps is extremal in the sense of having the smallest dilatation among all maps in its homotopy class;
\item
 these representatives are unique up to composition with conformal maps;
\item any topological self-mapping of a Riemann surface is homotopic to a map in this class.
\end{enumerate}

This is a formulation of what became later known as the \emph{Teichm\"uller theorem}.\index{Teichm\"uller theorem}\index{theorem!Teichm\"uller}

Teichm\"uller states that this result will also allow the computation of distances between points in the moduli space.

\section{Definition of principal regions through metrics}

 This   is the first technical section of the paper and, in fact, the paper remains technical until Section \ref{section23}. The idea in this section is to translate problems concerning Riemann surfaces into problems concerning Riemannian metrics. This is one of the most important sections in the paper, because Teichm\"uller sets here the tools for the results he will prove in the following sections.

In \S\,53, he explains how a sufficiently regular homeomorphism from a fixed surface $\mathfrak{h}_0$ onto another Riemann surface $\mathfrak{h}$  induces a Riemannian metric $ds^2$ on $\mathfrak{h}_0$ which is well defined up to multiplication by a scalar function. The Riemannian metric is written, in the tradition of Gauss, as
\[ds^2=\lambda(Edx^2+2Fdxdy+Gdy^2),
\]
with $\lambda$ denoting a scalar function.
The smoothness of the Riemannian metric depends on the smoothness of the homeomorphism which we started with, and it is not specified.

 Teichm\"uller shows how the Riemannian metric stems naturally from the conformal point of view. He starts with a fixed principal region $\mathfrak{h}_0$ and a homeomorphism onto another Riemann surface $\mathfrak{h}$. If $z=x+iy$ is a holomorphic local coordinate around a point in $\mathfrak{h}_0$ and  $z'=x'+iy'$ a holomorphic local coordinate around the image point in $\mathfrak{h}$, then, $z'$, as a function of $z$, in general, does not satisfy the Cauchy-Riemann equations. However, one can write:
\[\vert dz'\vert^2=dx'^2+dy'^2=Edx^2+2Fdxdy+Gdy^2,
\]
where $E,F,G$ are expressed as follows, in terms of the partial derivatives of $x'$ and $y'$ with respect to $x$ and $y$:
\[E=\left(\frac{\partial x'}{\partial x}\right)^2+
\left(\frac{\partial y'}{\partial x}\right)^2,\]
\[F=\frac{\partial x'}{\partial x}
 \frac{\partial x'}{\partial y}
 +\frac{\partial y'}{\partial x}
 \frac{\partial y'}{\partial y}
,\]
and
\[G=\left(\frac{\partial x'}{\partial y}\right)^2+
\left(\frac{\partial y'}{\partial y}\right)^2.\]
A class of metrics of the form 
\[ds^2=\lambda \vert dz'\vert^2=\lambda(Edx^2+2Fdxdy+Gdy^2),\]
with $\lambda$ a positive function,
is then naturally considered on $\mathfrak{h}_0$. This class of metrics is invariant under change of the local parameter $z$. Teichm\"uller recalls that conversely, the Riemann surface  $\mathfrak{h}$ together with its marking is uniquely determined by a metric $ds^2$ on $\mathfrak{h}_0$, which is well defined up to the scalar factor. In the smooth case, this result is due to Gauss, and it is equivalent to the existence of the so-called \emph{isothermal coordinates}.\index{isothermal coordinates} More precisely, Gauss' result says that any oriented surface equipped with a (smooth) Riemannian metric $g$ admits a unique complex structure with local parameters $z$ such that the metric can be written as $g(z)\vert dz\vert^2$, with $g(z)$ a positive smooth function defined in the local charts.

In the absence of smoothness, some regularity conditions are needed, and Teichm\"uller does not specify them. He will specify them in \cite{T29}.


The representation of conformal structures by Riemannian metrics will be very important later on in the paper, in particular, in the discussion of the infinitesimal theory of quasiconformal mappings.

In \S\,54, Teichm\"uller defines the equivalence relation between marked Riemann surfaces seen as equivalence classes of Riemannian metrics. This point of view is shown to be equivalent to the one of Riemanniann metrics on a fixed surface (the base surface). The equivalence relation leads to a definition of the space $R^\sigma$ as a space of equivalence classes of Riemannian metrics, where the equivalence is defined by pull-back by a map homotopic to the identity. The homeomorphisms that are used are termed ``sufficiently regular."

Teichm\"uller says that one difficulty concerning the topology of the spaces  $\mathfrak{R}^\sigma$ and $R^\sigma$ and  which was latent since \S\,13 is now solved since one can say that  two elements of such a space are close if corresponding metrics are close.

In \S\,55, he studies the effect of quasiconformal mappings on the metrics. He gives a new formulation of the existence problem for quasiconformal mappings, namely, finding the best quasiconformal mapping (an extremal one) homotopic to the identity, in terms of the deformations of the metrics. He sets the following new problem on the existence of extremal quasiconformal mappings:
\begin{problem}\label{probmetric}
Let two ``arbitrary"\footnote{Teichm\"uller puts the word ``arbitrary" in quotation marks.} metrics $ds_1^2$ and $ds_2^2$ on the principal region $\mathfrak{H_0}$ be given. Determine a mapping $A$ from $\mathfrak{H}_0$ onto itself which can be deformed to the identity such that the maximum of the dilatation quotient $D$ for the pair of metrics $ds_1^2$, $ds_2^2$ is smallest.
\end{problem}
Two equivalent ways of defining the space $R^\sigma$ show that this problem is equivalent to Problem \ref{prob3}.
Through Problem \ref{probmetric}, Teichm\"uller introduces the notion of dilatation quotient\index{dilatation quotient}\index{dilatation quotient!metric} relative to a pair of Riemannian metrics. He gives the definition in the two following subsections. One problem is now to compute the dilatation quotient in terms of the coefficients $E,F,G$ of the two metrics.

In \S\,56, he gives a formula for the dilatation quotient\index{dilatation quotient}  between two Riemann surfaces in terms of the coefficients of the Riemannian metrics and he notes that with this formula one can recover the generic formula for quasiconformal dilatation in terms of complex local parameters, established in \S\,15.

The formula is given by
$$
\Delta = \sqrt{EG-F^2};\, K=\frac{E+G}{2\Delta}; D=K+\sqrt{K^2-1}.
$$
This is found by a computation, using linear algebra.

In \S\,57, he summarizes the previous paragraphs and he gives the following precise form of the problem studied:
\begin{problem}
Let a principal region $\mathfrak{H}$ be given. One must find a set of metrics on it, of the form
$$
ds^2 = Edx^2 + 2F dxdy+Gdy^2 ,\;\left( z=x+\mathrm{i} y \textit{ is the local coordinate on } \mathfrak{H} \right)
$$
with the following three properties:
\begin{enumerate}
\item\label{item1} If one computes for an arbitrary metric $ds^2 =Edx^2 + 2F dxdy+Gdy^2$ a ``dilatation quotient" $D$ from the formulae
$$
\Delta = \sqrt{EG-F^2};\, K=\frac{E+G}{2\Delta}; D=K+\sqrt{K^2-1},
$$
then the metrics of the set have constant dilatation quotient, i.e. $D=$const.
\item\label{item2} When there exists, for a metric $ds^2$, a mapping $A$ from $\mathfrak{H}$ onto itself which can be deformed to the identity and for which $ds^2 \left( A\left( \mathfrak{p}\right)\right)$ belongs to this set and has constant dilatation quotient $K$, then the maximum of the dilatation quotients of $ds^2$ is itself
$\geq K$ with equality if and only if $ds^2 \left(A\left( \mathfrak{p}\right) \right)= ds^2 \left( \mathfrak{p}\right)$ holds.
\item\label{item3} To an ``arbitrary" metric $ds^2$ there always exists such an $A$.
\end{enumerate}
\end{problem}
Let us make some brief comments on this problem. Item (1) explains how to define, for a given Riemannian metric, the associated dilatation quotient.
 Item (2) implies that if there exists a map homotopic to the identity such that the dilatation quotient of its associated metric is constant, then the map is extremal. Item (3) says that there always exists an extremal map. Teichm\"uller adds that in the case of oriented closed surfaces, the quadratic differentials studied in Section \ref{sectionconj} provide examples of conformal metrics with constant dilatation.

In \S\,58, Teichm\"uller focuses on the link between a Riemannian metric and what he already introduced under the name ``direction field" (a singular foliation on the surface).\index{differential!quadratic} For this, he considers a Riemannian metric $ds^2$ on a principal region $\mathfrak{H}$,  and using complex coordinates $z$ on $\mathfrak{H}$   he rewrites (locally) this metric as
\begin{equation}\label{metriccompl}
ds^2 = \Lambda\vert dz \vert^2 +\mathrm{Re}\left( Hdz^2 \right),
\end{equation}
where $\Lambda$ is a positive function and $H$ a complex function, the two functions being determined by the coefficients $E$,  $F$ and $G$. From (\ref{metriccompl}), one can write the dilatation quotient in terms of $\Lambda$ and $\vert H \vert$. Teichm\"uller adds that the direction field is determined by the condition $Hdz^2 >0$. Thus, a metric $ds^2$ gives a dilatation quotient and a direction field and conversely, a dilatation quotient and a direction field determine a metric $ds^2$.

In \S\,59, Teichm\"uller studies the behaviour of $\Lambda$ and $H$ under a (conformal or anticonformal) change of coordinates.
First, he considers only a direct (conformal) coordinate change.
He finds that the expressions $\Lambda \vert dz\vert^2$ and $Hdz^2$ are invariant. (Note that the latter is a quadratic differential, but not necessarily holomorphic.)
He then introduces the quotient
\begin{equation}\label{beltrami1}
q=\frac{Hdz^2}{\Lambda \vert dz \vert^2},
\end{equation} which is invariant. Its conjugate $\overline{q}$ is what we call today a \emph{Beltrami differential}.\index{Beltrami differential} Furthermore, the dilatation quotient associated with the metric $ds^2$ is expressed in terms of $q$ by:
\begin{equation}\label{beltrami2}
D=\sqrt{\frac{1+\vert q \vert}{1-\vert q \vert}}.
\end{equation}
This equality is an analogue of the relation between the dilatation quotient and what is known today as the \textit{Beltrami coefficient}.\index{Beltrami coefficient} The last fact will be explained in the next section.

Teichm\"uller then makes similar considerations for anti-holomorphic coordinate changes.

In \S\,60, he gives explicitly, for a closed orientable surface,  the values of $H$, $\Lambda$ and $q$, in the case where the metric is given by a holomorphic (he says ``everywhere finite" ; cf. \S\,44) quadratic differential and a constant $K\geq 1$. He declares that this case is of particular interest. Indeed, we shall see that this will solve the problem addressed in \S\,57.

\section{Infinitesimal quasiconformal mappings}

In this section, Teichm\"uller introduces the notion of \textit{infinitesimal quasiconformal mapping}.\index{infinitesimal quasiconformal mapping} Such a mapping is given by a vector field on the surface. With this notion,  he establishes a new extremal problem, in \S\,66. This problem will lead him to solve an equation in Section \ref{chapitre19} (Equation (\ref{resolution}) below). This equation is a form of the \textit{Beltrami equation}.\index{Beltrami equation} He will also determine the tangent space of the Teichm\"uller space.

In \S\,61 and \S\,62, Teichm\"uller determines, with the previous notation, the expression of a metric after an infinitesimal perturbation. He fixes a metric $ds^2$  and he gives an expression for the parameters $q$ and $D$ of the metric $ds^2 + \delta ds^2$.  These parameters are given by a function $B$ which ``describes the infinitesimal quasiconformal mapping." The function $B$ alone determines the new metric obtained. As $B$ depends on local coordinates, Teichm\"uller rather considers the expression $B \frac{dz^2}{\vert dz \vert^2}$ which does not depend on coordinates. Eventually, $B$ will be seen to represent a tangent vector to Teichm\"uller space (in his notation, a point of $L^\sigma$).

In \S\,63, he studies how an infinitesimal mapping changes the function $B$, identifying a metric $ds^2(x)$ with $ds^2 \left( A(x)\right)$. He notes that the mapping $A$ is given in local coordinates by a complex function $w$ such that $\frac{w}{dz}$ does not depend on the change of coordinates. He adds some hypotheses on $w$ when he considers principal regions with non-empty boundary and/or with distinguished points. He writes:
\begin{quote}
 We thus obtain the most general infinitesimal mapping from $\mathfrak{H}$ onto itself with the help of a local function $w$ which still depends upon the parameter choice so that $\frac{w}{dz}$ vanishes at the distinguished inner and boundary points and is real along the boundary curves.
\end{quote}

In \S\,64, he determines the function which characterizes the infinitesimal metric obtained by an infinitesimal mapping $A$. He denotes this function by $B^*$ and he shows that
$$
B^* = B + \overline{w_x +\mathrm{i} w_y}.
$$
He then defines a new equivalence relation on metrics obtained by an infinitesimal perturbation of a given metric. Explicitly, he writes:
\begin{quote}
As a result we put in the same class all the $B\frac{dz^2}{\vert dz \vert^2}$ on $\mathfrak{H}$ that arise from each other through the modification
$$
B^* = B + \overline{w_x + \mathrm{i} w_y},
$$
where $\frac{w}{dz}$ must satisfy the conditions collected at the end of \S\,63. This yields the problem of finding in each of these classes a $B\frac{dz^2}{\vert dz\vert^2}$ with the smallest possible maximum for $\vert B\vert$.
\end{quote}
He notes that these families correspond bijectively to a small neighborhood of a point in the Teichm\"uller space $R^\sigma$. He adds that the space of classes of $B\frac{dz^2}{\vert dz \vert^2}$ forms a ``linear manifold" (that is, a vector space) of rank $\sigma$. This fact will be proved in Sections 20 and 24. This is the vector space he denotes by $L^\sigma$.

In \S\,65, he shows that the expression $\overline{\left( w_x + \mathrm{i} w_y\right)}\frac{dz^2}{\vert dz \vert^2}$ does not depend on coordinate changes.

In \S\,66, he returns to the space $L^\sigma$ and defines a norm on it, denoted by $\|\cdot\|$. He admits that this norm separates points; he will show this in Sections 20 (for the closed surface case) and 24 (for the general case). The most interesting fact in this section is the statement of a new problem, which could be considered as equivalent to Problem \ref{probmetric} but from an infinitesimal point of view. The problem, as he states it, is the following:
\begin{problem}\label{probbeltramimain}
Let a principal region $\mathfrak{H}$ be given. To find a set of invariant $B\frac{dz^2}{\vert dz \vert^2}$ on it  which has the following properties:
\begin{enumerate}
\item for all $B\frac{dz^2}{\vert dz \vert^2}$ in this set, $\vert B\vert$ is invariant;
\item if $B^* \frac{dz^2}{\vert dz \vert^2}\sim B\frac{dz^2}{\vert dz \vert^2}$ and $B\frac{dz^2}{\vert dz \vert^2}$ belongs to this set, then
$$
\max{\vert B^* \vert}\geq \vert B \vert,
$$
with equality only in the case where $B^* \frac{dz^2}{\vert dz \vert^2}=B\frac{dz^2}{\vert dz \vert^2}$;
\item to each invariant $B^* \frac{dz^2}{\vert dz \vert^2}$ there is a $B\frac{dz^2}{\vert dz \vert^2}$ in the set with $B^* \frac{dz^2}{\vert dz \vert^2}\sim B\frac{dz^2}{\vert dz \vert^2}$.
\end{enumerate}
\end{problem}
Through this new problem, Teichm\"uller introduces the notion of ``extremal infinitesimal quasiconformal mapping."\index{extremal infinitesimal quasiconformal mapping} At the end of this subsection, he announces that the solution of Problem \ref{probbeltramimain}, for closed surfaces, is determined by regular quadratic differentials. This will be proved in the following three sections.

\section{Extremal infinitesimal quasiconformal mappings}\label{etudeprob}

In this section, and until Section 20, the author considers only the case of closed oriented surfaces. He shows that the set
$$\left\lbrace c \frac{d\zeta^2}{\vert d\zeta \vert^2}, \; c>0\textrm{ and } d\zeta^2 \textrm{ is a regular quadratic differential}\right\rbrace$$
satisfies the first two properties of the statement of Problem \ref{probbeltramimain}. The third property will be shown in Sections 19  and 20.

In \S\,67, he considers the case where $B\frac{dz^2}{\vert dz\vert^2}\sim 0$. In this case, the equation to be solved is:
\begin{equation}\label{eqbeltrami}
\overline{B}=w_x + \mathrm{i} w_y,
\end{equation}
which he transforms into
\[\overline{B}_x -i\overline{B}_y =w_{xx} + \mathrm{i} w_{yy},
\]
which bears some formal analogy with the Poisson equation.\index{Poisson equation} This analogy will be mentioned again in \S\,76, but
Teichm\"uller says that the usual theory developed for the solution of the Poisson equation is not useful in the present context.

 The exact solutions of Equation (\ref{eqbeltrami}) will be the topic of Section \ref{chapitre19}.

In \S\,68, he introduces the following new definition.
\begin{definition}\label{definitionlocallyextremalsmallscales}
The differential form $B\frac{dz^2}{\vert dz \vert^2}$ is said to be \emph{locally extremal}\index{locally extremal} if $\vert B \vert$ is constant and at every point $\mathfrak{p}$ where $B\frac{dz^2}{\vert dz\vert^2}$ is defined, there is a simply connected uniformizing neighborhood $\mathfrak{U}$ which is mapped onto the plane through $z$, with the following properties: $B\left( z \right)$ is defined inside and on the boundary of $\mathfrak{U}$ and if $w$ vanishes outside $\mathfrak{U}$ and $\vert B+ \overline{w_x +\mathrm{i} w_y}\vert\leq\vert B\vert$ holds all over $\mathfrak{U}$, then $\vert B+ \overline{w_x +\mathrm{i} w_y}\vert=\vert B\vert$.
\end{definition}

 Definition \ref{definitionlocallyextremalsmallscales} is some kind of definition of \textit{locally extremal infinitesimal quasiconformal mapping}.\index{extremal infinitesimal quasiconformal mapping!locally}

In \S\,69, Teichm\"uller shows, using an idea introduced in \S\,68, that $\frac{dz^2}{\vert dz\vert^2}$ is locally extremal. He writes: ``In this proof as in the whole paper, we have renounced to be precise about assumptions."

In \S\,70, he proves the following proposition.
\begin{proposition}\label{proposition}
 If $f\left( z \right)$ is analytic and different from $0$ and $\infty$, then
 $$
 \frac{f\left( z \right)}{\vert f\left( z \right)\vert}\frac{dz^2}{\vert dz\vert^2}
 $$
 is locally extremal.
\end{proposition}
Through Proposition \ref{proposition}, we start to see the link with \emph{holomorphic} quadratic differentials. This link is more explicit in the following subsection. In \S\,71, he obtains
\begin{proposition}
If $f\left( z \right)$ is regular analytic, then $\frac{f\left( z \right)}{\vert f\left( z \right)\vert}\frac{dz^2}{\vert dz\vert^2}$ is locally extremal. In particular for a non-zero, everywhere finite quadratic differential $d\zeta^2$ and a constant $c>0$,
$$
c\frac{d\zeta^2}{\vert d\zeta \vert^2}
$$
is always locally extremal.
\end{proposition}
The property of being ``everywhere finite" allows a quadratic differential to have zeros or poles of order at most one. Teichm\"uller writes the following:
\begin{quote}
We can actually carry over the proof given in \S\,70 unchanged if the analytic function $f(z)$ has zeros.   On the other hand, for a first-order pole, difficulties arise in the partial integration and, for a pole of higher order, all integrals diverge even more.
\end{quote}

In \S\,72 and \S\,73, Teichm\"uller shows in ``a heuristic manner" that if $B\frac{dz^2}{\vert dz \vert^2}$ is locally extremal, then there exists an everywhere finite quadratic differential $d\zeta^2$ and $c>0$ such that
$$
B\frac{dz^2}{\vert dz \vert^2}=c\frac{d\zeta^2}{\vert d\zeta \vert^2}.
$$

In \S\,74 and \S\,75, he shows the following:
\begin{proposition}\label{propextre}
Let $d\zeta^2$ be an everywhere finite quadratic differential and $c>0$. Then $c\frac{d\zeta^2}{\vert d\zeta\vert^2}$ is locally extremal.
\end{proposition}
Using Proposition \ref{propextre}, he shows that the set
$$\left\lbrace c \frac{d\zeta^2}{\vert d\zeta \vert^2}, \; c>0\textrm{ and } d\zeta^2 \textrm{ is a regular quadratic differential}\right\rbrace$$
satisfies Properties (1) and (2) of Problem \ref{probbeltramimain}.

In \S\,76, the goal is to prove that the previous set also satisfies Property (3) of Problem \ref{probbeltramimain}. For this, due to the linear aspect of the space $L^\sigma$, Teichm\"uller wishes to give a  necessary condition for an element to be equivalent to $0$. He  states:
\begin{theorem}\label{Theo}If $B\frac{dz^2}{\vert dz\vert^2}$ is invariant and
$$
\int{\overline{B}\frac{\boxed{dz}}{dz^2}d\zeta^2}=0
$$
  for any everywhere finite quadratic differential $d\zeta^2$ on $\mathfrak{F}$, then there is an invariant $\frac{w}{dz}$ with
\begin{equation}\label{resolution}
B=\overline{w_x + \mathrm{i} w_y}.
\end{equation}
\end{theorem}

In Teichm\"uller's notation, $\boxed{dz}$ is an infinitesimal area form. In more conventional terms, $\boxed{dz}=\frac{1}{2\mathrm{i}}dz\wedge d\overline{z}$.

Teichm\"uller proves the converse of this theorem, in order to show Property \ref{propextre}. He says about Theorem \ref{Theo}:  ``I have originally taken this proposition as a postulate, then proceeded further and, much later, looked for a proof of it."
The proof of this theorem, which is the \textbf{key} of Teichm\"uller's theorem, will be established in the following section.

\section{The equality $\overline{w_x + \mathrm{i} w_y }= B$}\label{chapitre19}

As previously mentioned, the aim of this section is to show Theorem \ref{Theo}. Teichm\"uller divides this section into three parts in which he shows the theorem in the cases of the Riemann sphere, of the torus and of closed hyperbolic surfaces respectively.

From \S\,77 to \S\,80, he considers the Riemann sphere. There is no everywhere finite quadratic differential on this surface, but he shows a related result, namely, he produces a solution $w$ to the equation
$$
\overline{w_x + \mathrm{i} w_y}=B,
$$
assuming that the function $B$ on the sphere is sufficiently regular near $0$ and $\infty$. 
Moreover, he shows that $w$ is uniquely determined if we suppose that $w$ vanishes at $0$, $1$ and $\infty$.

Using his notation, the function $w$ is given by
\begin{equation}
w\left( z \right)=\frac{1}{2\pi}\int_{\mathbb{C}\mathrm{P}^1}{W\left( z, z_0 \right)\overline{B\left( z_0 \right)}\boxed{dz_0}}
\end{equation}
where the integral is taken over the sphere
and where $W\left( z, z_0 \right)$ is the kernel
$$
W\left( z, z_0 \right)= \frac{z\left( z-1 \right)}{z_0 \left( z_0 -1 \right) \left( z-z_0 \right)}.
$$
 
Moreover, the function $w(z)$ is unique, assuming that the inverse differential $\frac{w(z)}{dz}$ vanishes at $0,1,\infty$.

From \S\,81 to \S\,83, he gives a proof of Theorem \ref{Theo} in the case of the torus. He passes to the universal cover, which is the plane. If $u$ is a holomorphic coordinate in the plane, then an everywhere finite quadratic differential is of the form $adu^2$, with $a$ a complex number (since the complex dimension of the space $\mathfrak{M}^2$ is one), and an everywhere finite inverse differential is of the form $\frac{a}{du}$, with $a$ again a complex number (the complex dimension of the space $\frac{1}{\mathfrak{M}}$ is also one). He finds a solution which bears some similarities with the one of the sphere decribed above,  and where the integration is over a period parallelogram in the universal cover. 
Here, $B$ is an invariant differential, which is periodic with periods $w_1$ and $w_2$ associated with the torus, and such that
 \[\int\int B(u)\boxed{du}=0,
\]
the integration being again over a fundamental domain (a period parallelogram).

The solution $w$ of (\ref{resolution}) is then given as an integral of a kernel, as in the previous case, and it involves the \textit{Weierstrass $\zeta$-function}\index{Weierstrass $\zeta$-function} associated with the periods $w_1,w_2$. This is a rational function in $u$, with first order poles at $mw_1+nw_2$ ($m$ and $n$ integers), 
 and equivariant in $w_1$ and $w_2$, that is, satisfying 
 \[\zeta(u+w_i)=\zeta(u)+\eta_i
 \]
 for $i=1,2$.

This equivariance, the kernel, and the integral solution have their analogues in the case of closed oriented surfaces of genus $g\geq 2$, where the universal cover is the hyperbolic plane. This is treated in \S\,84 to \S\,87. Teichm\"uller proves the theorem using \textit{Poincar\'e series}.\index{Poincar\'e series} To be more precise, he considers the universal covering $\mathbb{D}$ of the surface $\mathfrak{F}$ and its associated \textit{Fuchsian group} $\mathfrak{G}$, he lifts the function $B$ to the disk (he denotes this lift again by $B$) and he considers the function
$$
w\left( \eta \right)=\frac{1}{2\pi}\int_{\mathbb{D}}{\frac{1}{\eta-\eta_0}\overline{B\left( \eta_0 \right)}\boxed{d\eta_0}},$$
the integration being again over a fundamental domain.
He shows that this function satisfies the condition to descend on the surface and satisfies Equation (\ref{resolution}).\footnote{This kind of approach is called ``coset enumeration" or ``telescoping sum corresponding to a  group action." It was later developed by Hejhal and Wolpert.}

We recall that at this point, Teichm\"uller does not solve Problem \ref{probbeltramimain}. He just poves Theorem \ref{Theo}.

\section{The linear metric space $L^\sigma$ of the classes of infinitesimal quasiconformal mappings}\label{vectorspace}

The main aim of this section is to finish the solution of Problem \ref{probbeltramimain} by using Theorem \ref{Theo}. We recall that if $B\frac{dz^2}{\vert dz \vert^2}$ is extremal in the sense of Definition \ref{definitionlocallyextremalsmallscales}, then there exists an everywhere finite quadratic differential $d\zeta^2$ and  $c>0$ such that $B\frac{dz^2}{\vert dz\vert^2}=c\frac{d\zeta^2}{\vert d\zeta \vert^2}$. To prove Property (3) of Problem \ref{probbeltramimain}, Teichm\"uller shows that the space $L^\sigma$ introduced in Section 17 is a vector space of real dimension $\sigma$. He also proves what he already announced in Section 17, namely, that $\| \cdot \|$ is a norm on $L^\sigma$.

In \S\,88, he gives a characterization of the class of $B\frac{dz^2}{\vert dz \vert^2}$. Indeed, Theorem \ref{Theo} implies that if $B^* \frac{dz^2}{\vert dz \vert^2}\sim B\frac{dz^2}{\vert dz \vert^2}$, then
\begin{equation}\label{pairing}
\forall 1\leq \mu\leq \tau,\;\int_{\mathfrak{F}}{\overline{B^* \left( z \right)}\frac{d\zeta_{\mu}^2}{dz^2}\boxed{dz}}=\int_{\mathfrak{F}}{\overline{B \left( z \right)}\frac{d\zeta_{\mu}^2}{dz^2}\boxed{dz}}
\end{equation}
where $\tau$ is the complex dimension of the space of everywhere finite quadratic differentials and $\left( d\zeta_{\mu}^2 \right)_{1\leq \mu \leq \tau}$ is a basis of this space. He considers a \textit{pairing} between Beltrami forms and quadratic differentials. Equality (\ref{pairing}) allows to associate, for each class of $B\frac{dz^2}{\vert dz \vert^2}$, $3g-3$ complex numbers $\left( k_{\mu} \right)_{1\leq \mu \leq \tau}$. He says that the natural question is whether the converse is true. The response is given by the following proposition which will be proved in the two subsequent subsections.

\begin{proposition}\label{propintermediaire}
Given any $\tau$ complex numbers $k_1 , k_2 , \cdots , k_\tau$, there always exists one and only one class of invariants $B\frac{dz^2}{\vert dz\vert^2}$.
\end{proposition}
In fact, Teichm\"uller will prove more. With this property, it is clear that $L^\sigma$ is a real vector space of dimension $6g-6$ (it is even, naturally, a complex vector space of dimension $3g-3$).

In \S\,89 and \S\,90, Teichm\"uller proves Proposition \ref{propintermediaire}. He shows that for each $\left( k_\mu \right)_{1\leq \mu \leq \tau}$ there exists a unique $c>0$ and a quadratic differential $d\zeta^2$ such that $c\frac{d\zeta^2}{\vert d\zeta \vert ^2}$ gives $\left( k_\mu \right)_{1\leq \mu \leq \tau}$. Thus, Property (3) of Problem \ref{probbeltramimain} is proved and therefore, the whole problem is solved.

In \S\,91, Teichm\"uller recalls briefly the definition of $\| \cdot \|$ on $L^\sigma$ and he proves that this is a norm. He also shows that for a vector $\mathfrak{t}=\left( k_\mu \right)_{1\leq \mu \leq 3g-3}$ and a complex number $a$, we have $\|a\cdot\mathfrak{t}\|=\vert a\vert \cdot\|\mathfrak{t}\|$. With this in mind, he asks whether this norm comes from a Hermitian structure. We quote his question:
\begin{quote}
[..]  I know about these linear metric spaces nothing more. 
In particular, it is completely uncertain whether there exists for instance a Hermitian matrix $(h_{\mu\nu})=(\overline{h_{\nu\mu}})$ such that
$$
\Vert\frak{k}\Vert{=}\sqrt{\sum_{\mu,\nu}\overline{k_{\mu}}h_{\mu\nu}k_{\nu}}.
$$
\end{quote}
 \section{Doubling}\label{sectionduplication}\index{doubling}

In this section, Teichm\"uller considers the cases of oriented bordered Riemann surfaces and non-orientable Riemann surfaces. He defines the associated doubled surfaces. These are closed Riemann surfaces with genus depending on the original surfaces. He introduces the notions of meromorphic functions, \textit{n-differentials} and \textit{divisors} on bordered non-orientable Riemann surfaces. These notions will be useful to establish the Riemann-Roch theorem for such surfaces, equipped with distinguished points. The theorem plays an important role in the next section, where the dimension formula announced in Section \ref{dimension}, Equation (\ref{eq1}), is proved.

In \S\,92, Teichm\"uller explains how to define the double\index{doubled surface} of a bordered or of a non-orientable Riemann surface.
For a bordered surface $\mathfrak{M}$, he defines the associated \textit{mirror surface}\index{mirror surface}  $\overline{\mathfrak{M}}$. This  is a bordered Riemann surface obtained by composing every local coordinate of $\mathfrak{M}$ by complex conjugation. To get the doubled surface, $\mathfrak{M}$ and $\overline{\mathfrak{M}}$ are glued along their boundary components. Teichm\"uller denotes this surface by $\mathfrak{F}$. We shall denote it here by $\mathfrak{M}^d$.

For a non-orientable surface $\mathfrak{M}$, its $2$-sheeted cover is likewise considered. It is a surface which has twice the number of boundary components than $\mathfrak{M}$.

In \S\,93, Teichm\"uller recalls the notion of \textit{algebraic genus}\index{algebraic genus} for bordered or non-orientable Riemann surfaces. This is the genus of the corresponding doubled surface. Let us recall  the formula.

 Let $\mathfrak{M}$ be a surface (bordered or non-orientable) whose topological data are $\left( g, \gamma, n \right)$, using the notation of Section \ref{dimension}. Then, $\mathfrak{M}^d$ is a closed surface whose genus is equal to
\begin{equation}\label{genrealgebrique}
G=2g+\gamma+n-1.
\end{equation}
This quantity $G$ is called the algebraic genus of $\mathfrak{M}$.

In \S\,94, Teichm\"uller defines divisors on $\mathfrak{M}$ and their associated degrees. He uses a particular normalisation in order to associate a divisor to $\mathfrak{M}^d$ with the same degree. To be more precise, he assumes that the integer associated with an inner point is even.

In \S\,95, he defines functions on $\mathfrak{M}$. He uses the term  \textit{function} instead of \textit{meromorphic function}. Meromorphic functions are real along boundary components. For bordered Riemann surfaces, he defines functions as follows:
\begin{quote}
Functions on the oriented bordered finite Riemann surface $\frak{M}$ are exactly the functions $f$ that are, up to poles, regular analytic in $\frak{M}$ and that are real on the boundary of $\frak{M}$.
\end{quote}
In the non-orientable case, the definition is:
\begin{quote}
 Functions $f$ defined on non-orientable finite surfaces $\frak{M}$ which depend on the choice of parameters so that $f$ is invariant under direct conformal parameter transformations but goes to $\bar{f}$ under indirect conformal transformations, which up to poles depend analytically on local parameters and which are real on the boundary, are exactly the functions on $\frak{M}$.
\end{quote}
Functions on $\mathfrak{M}$ correspond bijectively to \textit{symmetric} functions on $\mathfrak{M}^d$. They form a field on $\mathbb{R}$, or according to his terms, a \textit{real algebraic function field}.\index{real algebraic function field} This is in the tradition of Riemann, who identified a Riemann surface with the field of meromorphic function that it carries.

In \S\,96, Teichm\"uller defines for $n\in\mathbb{Z}$, \textit{$n$-differentials} on the general surface $\mathfrak{M}$ and for $n=-1$ (resp. $n=2$) he calls them inverse differentials (resp. quadratic differentials). Differentials have to be real on the boundary components. He recalls how an $n$-differential defines a divisor on $\mathfrak{M}$ and he also explains the bijective correspondance between differentials on $\mathfrak{M}$ and symmetric differentials on $\mathfrak{M}^d$. By his normalisation for a divisor on $\mathfrak{M}$, associated degrees are equal. As a consequence, the degree of any $n$-differential ($n\geq 0$) is equal to
$$
2 n \left( G-1 \right).
$$

In \S\,97, by the doubling process, he proves the Riemann-Roch theorem\index{Riemann-Roch theorem}\index{theorem!Riemann-Roch} for divisors on $\mathfrak{M}$. 

In \S\,98, he considers arbitrary topological types of surfaces (i.e with inner or boundary distinguished points) and proves, using the notation of Section \ref{dimension}, the following version of the Riemann-Roch theorem which is analogous to the dimension formula he gives in \S\,14:
\begin{quote}
 The difference between the maximal number of real-linearly independent quadratic differentials that have at most first-order poles at the distinguished points and the maximal number of real-linearly independent everywhere finite inverse differentials that vanish at the distinguished points is always equal to
\begin{equation}\label{formuledimension}
-6+6g+3\gamma+2h+k.
\end{equation}

\end{quote}
The next three sections contain a proof of the dimension formula and solve Problem \ref{probbeltramimain} in the general case.

\section{Regular quadratic differentials and inverse differentials of a principal region}

This section follows the same idea as Section \ref{etudeprob}. The author justifies the link between $\rho$ and the real dimension of everywhere finite reciprocal differentials. He  shows that the space of quadratic differentials satisfies the first two properties of Problem \ref{probbeltramimain}.

In \S\,99, he considers the vector space of \textit{regular inverse differentials},\index{regular inverse differential} that is, inverse differentials vanishing at the distinguished points. He establishes the link between the real dimension of this vector space and the value $\rho$.

In \S\,100, he recalls when $B\frac{dz^2}{\vert dz \vert^2}$ and $B^*\frac{dz^2}{\vert dz \vert^2}$ are equivalent and he says that the main aim is to find in each class the one with the smallest norm. He solves Problem \ref{probbeltramimain} in some particular cases   (doubly connected domains, quadrilaterals, sphere with four distinguished points). The solutions are given by quadratic differentials with conditions on the boundary and at distinguished points. These conditions lead to the definition of a general  \textit{regular quadratic differential},\index{regular quadratic differential} which is an everywhere finite quadratic differential ``apart from possible first order poles at the distinguished points." He  states the following theorem:
\begin{theorem}
The set
\[\left\lbrace c \frac{d\zeta^2}{\vert d\zeta \vert^2}, \; c\geq0\textrm{ and } d\zeta^2 \textrm{is a regular quadratic differential}\right\rbrace
\]
 solves Problem \ref{probbeltramimain}.
\end{theorem}
This theorem will be completely proved in Section 24. In \S\,101, Teichm\"uller proves, with the same method as in \S\,74, that the previous set satisfies the first two properties of Problem \ref{probbeltramimain}.

\section{The equation $\overline{w_x +\mathrm{i} w_y}=B$ for arbitrary principal region}

This section is devoted to the proof of a theorem analogous to Theorem \ref{Theo}. This theorem leads to the solution of the third property of Problem \ref{probbeltramimain} in the case of bordered or non-orientable surfaces with distinguished points. The idea is simple. The author starts by proving the theorem (Theorem \ref{Theo2} below) in the case of closed surfaces with distinguished points. After that, to treat the most general surfaces, he considers the corresponding doubles. By \textit{symmetrization}, he gets the general solution.

In \S\,102, Teichm\"uller establishes the following analogue of Theorem \ref{Theo}:
\begin{theorem}\label{Theo2}
If $B\frac{dz^2}{\vert dz\vert^2}$ is invariant up to complex conjugation and if
\begin{equation}\label{hypothese}
\int_{\mathfrak{M}}{\mathrm{Re}\left( \overline{B}\frac{d\zeta^2}{dz^2} \right)\boxed{dz}}=0
\end{equation}
for all regular quadratic differentials $d\zeta^2$ of the principal region $\mathfrak{H}$, then there exists a regular inverse differential $\frac{w}{dz}$ such that
\begin{equation}\label{eqbeltramigeneral}
B=\overline{w_x +\mathrm{i} w_y}.
\end{equation}
\end{theorem}
The theorem was already proved for closed oriented surfaces without distinguished points. The author provides arguments for the proof in the case of the punctured torus, and for a sphere with one, two or three distinguished points. He explains that we have to solve the corresponding equation in the absence of distinguished points and then use the general form of the solution, which is unique up to constants. Furthermore, these constants may be chosen in order for the solution $w$ to vanish at the distinguished points. Therefore, he gets a regular inverse differential.

In \S\,103, he proves Theorem \ref{Theo2} for the sphere with four distinguished points. For this, he uses a two-sheeted covering space (i.e a torus) ramified at four distinguished points and he obtains an inverse differential on it. He modifies this differential so that it descends to the sphere with four distinguished points and he checks that it vanishes at these points.

In \S\,104, he considers the most general case of a hyperbolic closed surface with distinguished points. He lifts the problem to the universal covering, and he recalls the form of the solution $w$. This solution uses Poincar\'e series. He proves that it vanishes at the distinguished points.
At this step, he proves Theorem \ref{Theo2} for closed surfaces with distinguished points.

In \S\,106, he considers bordered and non-orientable surfaces. Let $\mathfrak{M}$ be such a surface with $h$ inner distinguished points and $k$ boundary points. Thus, the corresponding double $\mathfrak{M}^d$ has $2h+k$ distinguished points. Teichm\"uller transfers the problem on $\mathfrak{M}$ to a problem on the double. The invariant $B\frac{dz^2}{\vert dz\vert^2}$ gives a symmetric invariant function on the double. By Equation (\ref{hypothese}), Teichm\"uller gets all the conditions to find a solution $w$. He modifies this solution in order to obtain a symmetric regular inverse differential and, therefore, a solution to Equation (\ref{eqbeltramigeneral}). Theorem \ref{Theo2} is now proved.

In \S\,107, he gives explicitly the solution of (\ref{eqbeltramigeneral}) in the case of a simply connected domain. He uses a solution already given in \S\,80.

\section{The dimension of the linear metric space $L^\sigma$}

In this section, Teichm\"uller finishes the solution of Problem \ref{probbeltramimain}. For this need, he proves that $\left(L^\sigma , \| \cdot \|\right)$ is a real normed vector space of dimension $\sigma$. The poof follows the scheme of Section \ref{vectorspace}.

In \S\,108, he remarks that, according to Theorem \ref{Theo2} (\S\,88), a class of invariants $B\frac{dz^2}{\vert dz \vert^2}$ is characterized by a vector $\left( k_\mu \right)_{1\leq \mu\leq \sigma} \in \mathbb{R}^{\sigma}$. Unlike  in \S\,88, he sets, for an invariant $B\frac{dz^2}{\vert dz\vert^2}$, and  $1\leq \mu\leq \sigma$,
\begin{equation}\label{vecteur}
k_\mu =\int_{\mathfrak{M}}{\mathrm{Re}\left( \overline{B}\frac{d\zeta_\mu^2}{dz^2} \right)\boxed{dz}}, 
\end{equation}
where $\left( d\zeta_\mu\right)_{1\leq \mu\leq \sigma}$ is a basis of the real vector space of regular quadratic differentials. Thus, he shows that $L^\sigma$ is a real vector space. In \S\,109, he proves the following:
\begin{proposition}
Let $\mathfrak{t}\in\mathbb{R}^\sigma -\left\lbrace 0\right\rbrace$. Then, there exists a unique $c>0$ and a unique regular quadratic differential $d\zeta^2$ of norm $1$ such that  $c\cdot d\zeta^2$ determines (by Formula (\ref{vecteur})) the vector $\mathfrak{t}$.
\end{proposition}
He notes that with this proposition, Property (3) of Problem \ref{probbeltramimain} is proved and therefore this problem is solved. We also see that the space $L^\sigma$ is a real vector space of dimension $\sigma$. He also checks that the norm $\| \cdot \|$ is indeed a norm.

In \S\,110, he shows that $L^\sigma$ is a complex vector space when $\mathfrak{M}$ is a closed oriented surface with $h$ distinguished points. 

In \S\,111, he finishes the proof of the dimension formula
$$
\sigma - \rho = -6+6g+3\gamma+3 n +2h +k.
$$
It remains to prove that Teichm\"uller space is a manifold of dimension $\sigma$. This will be done in the following section.

In \S\,112, the author gives a heuristic way to prove the dimension formula. Following an idea of Klein, he establishes this formula by induction on the number of distinguished points.

\section{Passing to finite mappings. $R^\sigma$ as a Finsler space}\label{section23}

This section may be seen as the most important section of this article. Teichm\"uller establishes  what we call today Teichm\"uller's theorems. For this, he makes the link between extremal infinitesimal quasiconformal mappings and extremal quasiconformal mappings and he shows that the Teichm\"uller space with respect to the Teichm\"uller metric is a Finsler space. He also constructs isometric embeddings of the unit disk into Teichm\"uller space.

In \S\,113, he considers an extremal mapping $f: \mathfrak{H }\rightarrow \mathfrak{H}^\prime$ where $\mathfrak{H}$ and $\mathfrak{H}^\prime$ are 
principal regions and a metric $ds^2$ on $\mathfrak{H }$, associated with $f$, written in complex coordinates as 
$$
ds^2 =\Lambda \vert dz \vert ^2+\mathrm{Re}\left( Hdz^2 \right).
$$
We recall that the direction field of $f$ is given by the condition $Hdz^2 >0$ and its quasiconformal dilatation is $\sqrt{\frac{\Lambda + \vert H\vert}{\Lambda-\vert H \vert}}$. 

Teichm\"uller considers the family of metrics
$$
ds_{t}^2 =\Lambda \vert dz \vert^2 +t\cdot \mathrm{Re}\left( Hdz^2 \right) \;\;\;\left( 0\leq t \leq 1 \right);
$$
whose associated quasiconformal mappings, denoted by $f_t$, are extremal. (This is justified in \S\,39.) He notes that for $t\sim 0$, $f_t$ becomes an infinitesimal extremal quasiconformal mappings. Thus, by the previous section, there exists $c>0$ and a regular quadratic differential $d\zeta^2$ such that
\begin{equation}\label{eq:qc}
ds^2 = \lambda\left( \vert d\zeta^2 \vert+ c \mathrm{Re}\left( d\zeta^2 \right) \right).
\end{equation}
The conjecture on quasiconformal mappings is now equivalent to the fact that the extremal (finite) quasiconformal mappings are always given by metrics of the form  (\ref{eq:qc}) where $d\zeta^2$ is a quadratic differential, $\lambda$ a positive function and $c$ a positive constant. Concerning this passage, Teichm\"uller says: ``\textit{This conjecture is based on nothing.}"\footnote{The equivalence between (uniquely) extremal mappings in a class and  (uniquely) extremal in the corresponding infinitesimal class, in the general case, was obtained later on by works of Hamilton, Krushkal', Reich, and Strebel.}

In \S\,114, he formulates his previous observation as the following theorem.  \begin{theorem}\label{existence}\index{Teichm\"uller existence theorem}\index{theorem!Teichm\"uller existence}
Every extremal quasiconformal mapping from some topologically fixed principal region $\mathfrak{H}$ to some principal region $\mathfrak{H}^\prime$ of the same fixed topological type corresponds to a metric of the form
\begin{equation}\label{metriqueextremale}
ds^2 = \lambda\left( \vert d\zeta^2 \vert + c \mathrm{Re}\left( d\zeta^2 \right) \right).
\end{equation}
Here $d\zeta^2$ is a regular quadratic differential on the principal region $\mathfrak{H}$, $c$ is a positive constant which must incidentally be taken $<1$ and $\lambda$ is a variable positive factor.
\end{theorem}
This theorem is a form of the Teichm\"uller existence theorem.\index{Teichm\"uller existence theorem}\index{theorem!Teichm\"uller existence} By a clever rewriting of (\ref{metriqueextremale}), he observes that the corresponding metric determines (up to a conformal mapping) a mapping in the coordinates given by $d\zeta^2$:
\begin{equation}\label{teichmullermap}
\zeta \mapsto K\mathrm{Re}\left( \zeta \right) + \mathrm{i} \mathrm{Im}\left( \zeta\right).
\end{equation}
In other words, the metric determines the Teichm\"uller mapping.

In \S\,115, Teichm\"uller declares that the mappings in (\ref{teichmullermap}) are extremal. This will be proved in Section 27. These mappings  solve partially Problem  \ref{probexistence&unicite} and also the equivalent Problem \ref{probmetric}. At this step, it solves it ``partially" because the author does not give any proof of the fact that ``every topological mapping from a principal region $\mathfrak{H}$ onto some other principal region $\mathfrak{H}^\prime$ can be deformed into a mapping" of a form given in (\ref{teichmullermap}). He says that he will give a proof, but it is not clear whether a proof is contained in the present paper.

In \S\,116, Teichm\"uller introduces a new metric which is a \textit{length metric}. For this he makes use in an essential way of the structure of the vector space $L^\sigma$. He denotes this new metric by $\left\lbrace \cdot , \cdot \right\rbrace$ and he explains why this metric induces on $R^\sigma$ the structure of a Finsler space.
In \S\,117, he proves:
\begin{theorem}
Let $P$ and $Q$ be two points in $R^\sigma$. Then
$$
\left\lbrace P, Q \right\rbrace= \left[ P Q \right].
$$
\end{theorem}

In \S\,119, he defines rays and shows that they are geodesics. Today, we call these rays \emph{Teichm\"uller rays}.\index{Teichm\"uller ray} He also says that from the point $P$, there \textit{emanates exactly one geodesic in every direction of $R^\sigma$}.

In \S\,120, he proves
\begin{theorem}
The geodesic rays corresponding to $d\zeta^2$ and $-d\zeta^2$ and that emanate from a point $P$ coalesce into one single geodesic line.
\end{theorem}

In \S\,121, he considers the space $R^\sigma$ for closed oriented surfaces with possibly distinguished points.
 Given an element $P$ of $R^\sigma$, a complex number $\varphi$ and $K\geq 1$, he considers the element $P(K,\varphi)$ of $R^\sigma$ whose direction is $e^{-i\varphi}d\zeta^2$ and at distance $\log K$ from $P$.
 The target principal region is denoted by $\mathfrak{h}(K,\varphi)$. It is describd by the metric
 \[
 \vert d\zeta^2\vert + \frac{K^2-1}{K^2+1}\mathfrak{R}e^{-i\varphi}d\zeta^2.\]
 He considers $\log K$ and $\varphi$ as ``geodesic polar coordinates" and
 the surface defined by all these $P(K,\varphi$) as a \emph{complex geodesic}.\index{complex geodesic} Such a surface is called today a \emph{Teichm\"uller disc}.\index{Teichm\"uller disc}
  In the special case where the surface $\mathfrak{h}$ is a marked torus with no distinguished point (\S\,26 to 28), the whole space $R^2$ is a ``unique complex geodesic." Equipped with the Teichm\"uler metric, this space $R^2$ is the hyperbolic plane (of constant curvature $-1$). More generally, Teichm\"uller establishes the following:
 \begin{theorem}
 Any complex geodesic is isometric to the hyperbolic plane.
 \end{theorem}

In \S\,122, he proves that $R^\sigma$ with respect to $[\cdot , \cdot]$ is uniquely geodesic. He uses for this the Teichm\"uller existence theorem (which he did not completely  prove yet).

In \S\,123, he uses polar coordinates in order to justify:
\begin{theorem}[Existence and uniqueness of geodesics between any two points]
 Between any two points of $R^\sigma$ passes a unique geodesic, and this geodesic is distance minimizing (not only locally).  
 \end{theorem}
This implies:
\begin{theorem}\label{theoremeteichmuller}
The space $R^\sigma$ of classes of conformally equivalent topologically fixed principal regions of a given type is homeomorphic to a $\sigma$-dimensional Euclidean space.
\end{theorem}
This theorem is related to the assertion made in \S\,13. It is probable that the reason of the use of $\mathfrak{R}^\sigma$ there, instead of $R^\sigma$, was to avoid introducing more definitions and notation. Theorem \ref{theoremeteichmuller} is sometimes called the \textit{Teichm\"uller theorem}.\index{theorem!Teichm\"uller}\index{Teichm\"uller theorem}

  \section{Simple cases}

In \S\,124 the author considers special cases of Theorem \ref{theoremeteichmuller}.  These cases are classified according to Klein's notion of algebraic genus\index{algebraic genus}\index{genus!algebraic} and reduced dimension\index{reduced dimension}\index{dimension!reduced} which were discussed in \S\,14. In \S\,124, he recalls Klein's classification of surfaces according to their algebraic genus. For each nonpositive integer $G$, he gives the number of bordered/non-orientable surfaces having algebraic genus $G$. (For closed surfaces, the number is 1.)

  There are 12 cases corresponding to the value $\tau=1$,  and he shows explicitly that in these cases $\sigma$ has the announced value. These are cases where the quadratic differential $d\zeta^2$ can be given explicitly.

 In \S\,128, he considers 11 special cases where $\tau=2$. He considers the torus as a two-fold  covering of the sphere ramified over four points. He declares: ``The structure of the linear space $L^\sigma$ is still not known to me for $\tau>1$."

\begin{figure}[ht!]
\centering
 \psfrag{0}{\small $0$}
 \psfrag{a}{\small $a$}
 \psfrag{b}{\small $a+ib_1$}
 \psfrag{c}{\small $a_1+ib_1$}
 \psfrag{d}{\small $a_1+ib$}
 \psfrag{e}{\small $ib$}
\includegraphics[width=.30\linewidth]{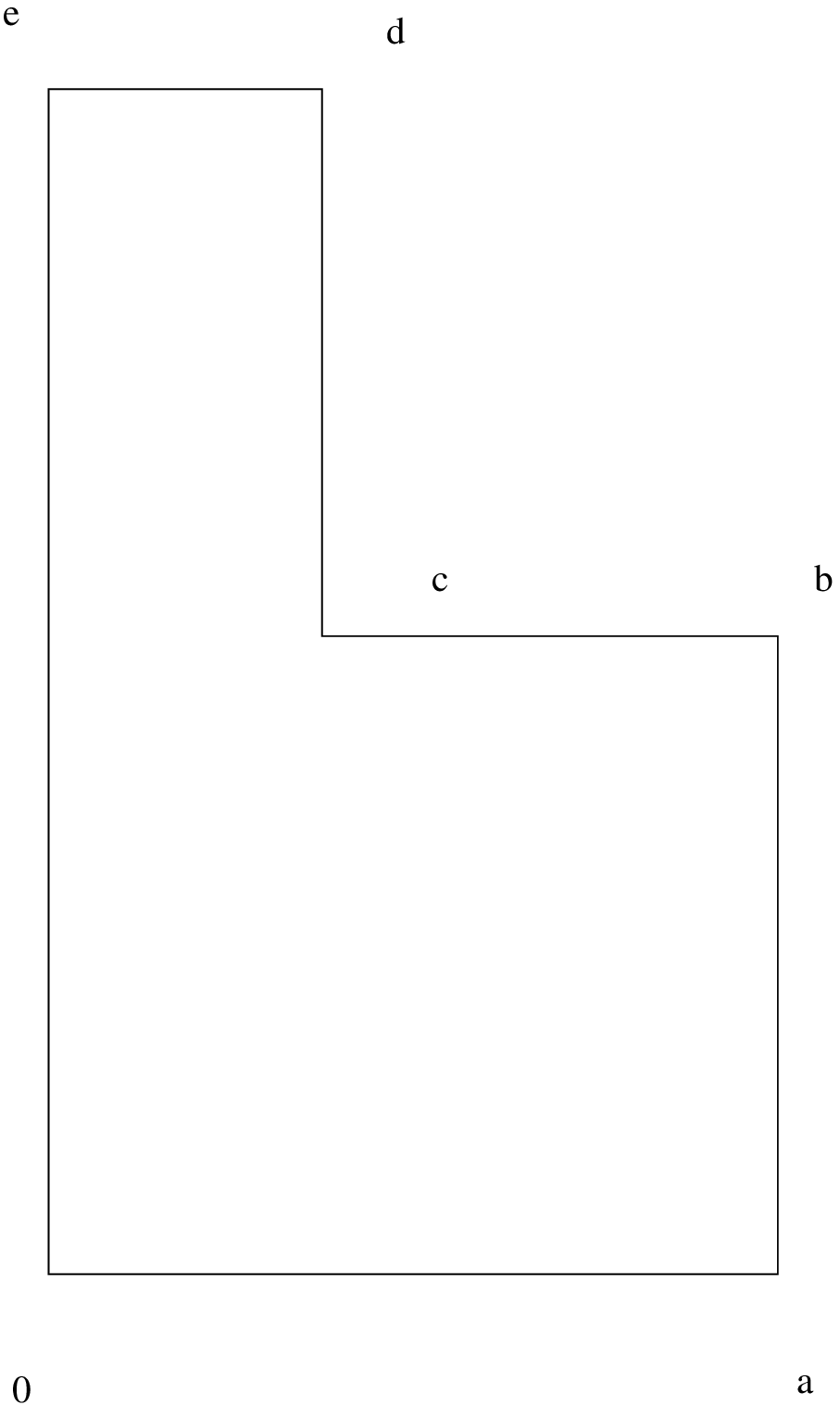}
\caption{\small{ {Hexagon}}}
\label{hexagon}
\end{figure}

 In \S\,129 he considers the extremal  quasiconformal mappings between pentagons, that is, principal regions with five distinguished points on the boundary. In this case, he can take explicit parameters by considering as a pentagon the upper half-plane with distinguished points $r_1,r_2,r_3,r_4,\infty$ and mapping it to the $L$-shaped Figure \ref{hexagon}. By convention, the vertices are only the ones with  acute angles (equal to $\pi/2$). 
 The space of regular quadratic differentials on such a principal region is of real dimension two, with  
possibly simple poles at the distinguished points. 
On such a pentagon, a map
 \begin{equation}\label{zeta}
 \zeta'=K\zeta+i\eta
 \end{equation}
 can be applied. Mapping again to the upper half-plane, the  map induced by (\ref{zeta}) is claimed to be extremal quasiconformal. Teichm\"uller treats completely the case of a pentagon in \cite{T24}.

 At the end of \S\,129, Teichm\"uller writes: ``It is advisable to carry out simpler examples of the same type before giving the general proof (in \S\S\,132-140)."

In \S\,130, he compares the dimension of the space of hexagons (with sides parallel to the axis
and with five distinguished pointing-out vertices) and the dimension of the Teichm\"uller space 
of  a principal region  given by
the upper half-plane with five distinguished points on the real axis. 
 Figure \ref{hexagon} is again used. The space of hexagons is three-dimensional. Two one-parameter families of hexagons inside this parameter space are highlighted.  One of them consists of Teichm\"uller geodesics. The author discusses the behavior at infinity of these geodesics. 
He also interprets the space of hexagons as a fiber bundle such that each fiber corresponds to a 
Teichm\"uller geodesic and each point in the base space corresponds to the two ends of a Teichm\"uller 
geodesic.

 Two spaces $R^2$ and $S^2$ are obtained by identification of points on the same curve. There is a \emph{contact transformation}\index{contact transformation} between the two spaces $R^2$ and $S^2$ (a transformation carrying lines to lines). Then follows a discussion of convergence of points to infinity.

   In \S\,131, he studies ring domains (annuli) with one distinguished boundary point on each boundary curve. These surfaces were already considered in \S\S\,3 and 33. An application of the Schwarz reflection principle leads to the torus.

\section{Proof of the extremality property}

In this section, Teichm\"uller gives a proof of two important results, namely,
\begin{itemize}
\item 
that Teichm\"uller maps are extremal in their homotopy classes;
\item that if a map, in the homotopy class of a Teichm\"uller map, has the same quasiconformal dilatation, then it is equal to the Teichm\"uller map up to composition by a conformal map.
\end{itemize}
The proof is long, and takes from \S\,132 to \S\,140. It is inspired by the case of the torus (given in \S\,28) but the technical details in the general case are much more involved. The proof uses geometry, topology and analysis. It makes use of the metric on the universal cover induced by the quadratic differential defining the Teichm\"uller map, the Hopf-Rinow theorem applied to the geodesics of that metric, the Gauss-Bonnet theorem for polygons with piecewise geodesic boundary, and other important tools.

In \S\,135, Teichm\"uller proves the following rigidity result  for the action of a conformal mapping on quadratic differentials: 
\begin{proposition}
A conformal mapping of a closed orientable Riemann surface which can be deformed to the identity can only multiply every quadratic differential by a constant of norm one.
\end{proposition}
 In \S\,136, he considers surfaces with boundary.

The results proved here are the so-called Teichm\"uller's uniqueness theorem.\index{Teichm\"uller uniqueness theorem}\index{theorem!Teichm\"uller uniqueness} Teichm\"uller's proof has been reproduced with more details by several authors, including Ahlfors and Bers.\footnote{Ahlfors, in his comments on his paper \emph{On quasi-conformal mappings} \cite{Ahlfors1954} in his \emph{Collected Works} \cite{Ahlfors:Collected} (Vol. 1, p. 1), writes: ``A complete proof of the uniqueness part of Teichm\"uller's theorem was included. Like all other known proofs of the uniqueness it was modeled on Teichm\"uller's own proof, which used uniformization and the length-area method.\index{length-area method} Where Teichm\"uller was sketchy I tried to be more precise." Bers, in his paper \cite{Bers1960} also gives a proof of uniquess, and he writes about it (p. 108): ``We merely re-word Teichm\"uller's own argument, making sure that it applies also to the more general definition of quasiconformality used here."}

\section{Conformal mappings of a principal region onto itself}

In \S\,141, Teichm\"uller proves the following:
\begin{proposition}
A conformal self-mapping of a closed orientable Riemann surface which can be deformed to the identity can be deformed to the identity through conformal maps.
\end{proposition}
In particular, for genus $>1$, since there are no holomorphic continuous deformations, the map itself is the identity. In other words he proves the following:
\begin{proposition} A conformal self-mapping of a closed orientable Riemann surface of genus $>1$ which is homotopic to the identity is the identity.
\end{proposition}

The proof uses the action on quadratic differentials.\footnote{Teichm\"uller's result was extended to an arbitrary quasiconformal self-mapping, see \cite{EMM}.}
Teichm\"uller's proof is based on the fact that when  $\mathfrak{H}$ is a closed surface of genus $g\geq 2$, the conformal mapping acts nontrivially on
the homology group $H_1(\mathfrak{H},\mathbb{Z})$ (or equivalently, on the space of holomorphic $1$-forms).
\footnote{Such an argument dates back to Hurwitz \cite{Hur}.}

In \S\,142, he introduces the \emph{mapping class group}, denoted by $\mathfrak{F}$, of a principal region $\mathfrak{H}$.
He explains that the definition of $\mathfrak{F}$ only depends on the topology of $\mathfrak{H}$,
and the result in \S\,141 implies that the quotient group $\mathfrak{R} | \mathfrak{C}$  (where
 $\mathfrak{R}$ is the conformal automorphism group of $\mathfrak{H}$ and $\mathfrak{C}$
 is the subgroup of  $\mathfrak{R}$ of conformal mappings homotopic to the identity)
can be seen as a subgroup of $\mathfrak{F}$.  Note that $\mathfrak{R} | \mathfrak{C}$ is finite, by the famous Hurwitz Theorem.\index{theorem!Hurwitz}

In \S\,143, Teichm\"uller considers the action of the mapping class group $\mathfrak{F}$ on the Teichm\"uller space
$R^{\sigma}$ and shows that this action preserves the Teichm\"uller metric.
He also defines the quotient group $\mathfrak{F} | \mathfrak{N}'$, where $\mathfrak{N}'$ is the group of
elements in $\mathfrak{F}$ that act trivially on $R^\sigma$. The quotient of $R^{\sigma}$ under $\mathfrak{F}$
is the moduli space $\mathfrak{R}^\sigma$. Teichm\"uller states that $\mathfrak{F} | \mathfrak{N}'$ acts on
$R^{\sigma}$ properly discontinuously.

In \S\,144, he shows that when $\mathfrak{H}$ is a torus,
$\mathfrak{F} | \mathfrak{R}'$ is identified with the classical modular group $\mathrm{PSL}(2,\mathbb{Z})$.
He points out that $\mathfrak{F} | \mathfrak{N}'$ contains non-trivial elliptic elements.

In \S\,145, Teichm\"uller considers closed surfaces of genus two.
In this case, each principal region $\mathfrak{H}$ is hyperelliptic and
 $\mathfrak{H}$ is represented by a two-sheeted
branched covering of  the Riemann sphere.
From here, he shows that the extremal quasiconformal mapping between any two principal regions
commutes with their involutions.
This implies that $ \mathfrak{N}'$ is isomorphic to $\mathbb{Z}_2$.

Given a pair of principal regions $\mathfrak{H}$ and $\mathfrak{H}'$, let $K$ (rep. $K'$) be a conformal automorphism of $\mathfrak{H}$ (rep. $\mathfrak{H}'$)
such that $K$ and $K'$ agree as elements of the mapping class group. Denote the extremal quasiconformal mapping
between $\mathfrak{H}$ and $\mathfrak{H}'$ by $E$ (here it  is assumed that $E$ arises from a
regular quadratic differential and from a dilatation constant).
In \S\,146, Teichm\"uller proves that the initial regular quadratic differential of $E$ is invariant under $K$.
Conversely, he shows that any regular quadratic differential on $\mathfrak{H}$ that is invariant under $K$
arises in this way. As a corollary, he obtains the interesting result that if an element of the mapping class group leaves two distinct points
in $R^\sigma$ fixed, then it fixes the Teichm\"uller geodesic (even the complex Teichm\"uller geodesic) through these two points.

In \S\,147, Teichm\"uller introduces the notion of \emph{geodesic manifold in $R^\sigma$}.\index{Teichm\"uller space!convexity}\index{convexity!Teichm\"uller space} In modern terms, this is a totally geodesic subset of $R^\sigma$, endowed with the Teichm\"uller metric.
It follows that the subset of $R^\sigma$ fixed by an element of $\mathfrak{F}$ (if not empty)
is a geodesic manifold, and the intersection of geodesic manifolds is again a geodesic manifold.
The dimension of the geodesic manifold fixed by a subgroup $\mathfrak{U}$ of $\mathfrak{R} | \mathfrak{C}$
is equal to maximal rank of linearly independent $\mathfrak{U}$-invariant regular quadratic differentials.
Teichm\"uller asks the question of computing the dimension algebraically.

In \S\,148, Teichm\"uller asks whether any finite subgroup $\mathfrak{U}$ of $\mathfrak{F}$ can be realized as a subgroup
of $\mathfrak{R} | \mathfrak{C}$ which acts on some principal domain $\mathfrak{H}$ conformally.
This is the famous Nielsen Realization Problem,\index{Nielsen Realization Problem}\index{problem!Nielsen realization} which was solved by Kerckhoff in 1983 \cite{kerckhoff1}.

In \S\,149, Teichm\"uller uses quasiconformal mappings to prove that given a topological self-mapping $S$ of  $\mathfrak{H}$
such that $S^2$ is homotopic to the identity, there exists a topological self-mapping
 of order two in the homotopy class of $S$. Then he explains that the Nielsen Realization Problem\index{Nielsen Realization Problem}\index{problem!Nielsen realization} for 
 a finite subgroup $\mathfrak{U}$ of $\mathfrak{F}$ is equivalent the problem of whether
 $\mathfrak{U}$ has a fixed point in Teichm\"uller space.

From \S\,150 to \S\,152, Teichm\"uller illustrates the preceding general discussion by considering a hyperelliptic surface of genus $g>1$ and a topological self-mappings $G$  satisfying $G^2=1$ ($G$ is called a \emph{distinguished mapping}).
In \S\,150, he shows that $G$ can be realized by some $F$ of order two,
which is conformal on some hyperelliptic principal region $\mathfrak{H}$ ($\mathfrak{H}$
is a two-sheet covering of the Riemann sphere, with $2g + 2$ ramification points, and $F$ acts on $\mathfrak{H}$
by exchanging the sheets). Then he computes the dimension of the geodesic manifold fixed under $F$. (Here $F$ is considered as an element of the mapping class group $\mathfrak{F}$.)

In \S\,151, he shows that any two homotopic distinguished mappings are conjugate
by a topological mapping homotopic to the identity.

In \S\,152, Teichm\"uller discusses the number of conjugacy classes of
$G$ in $\mathfrak{F}$.  As shown in \S\,150, this number is equal to the number of distinguished classes,
and it is also equal to the index of $[\mathfrak{F}, \mathfrak{N}_G]$, where $\mathfrak{N}_G$
is the normalizer of $G$ in $\mathfrak{F}$. He gives a geometric interpretation of $\mathfrak{N}_G$
in terms of the branched cover over the Riemann sphere. At the end of \S\,152, he proposes the following
three related topological problems:

\begin{enumerate} 
  \item Prove the result of \S\,151 in a purely topological way.
  \item Prove in a purely topological way that in the case $g = 2$, $[\mathfrak{F}, \mathfrak{N}_G]=1$ (as shown in \S\,145).
  \item Compute the index $[\mathfrak{F}, \mathfrak{N}_G]$ for general $g>1$.
\end{enumerate}

In \S\,153, Teichm\"uller considers branched coverings between finite Riemann surfaces. He
defines the notion of \emph{normal covering}\footnote{This is also known as Galois coverings.} between finite Riemann surfaces.
Then he mentions that a covering between finite Riemann surfaces,
say  $f: \mathfrak{M}\to \mathfrak{M}^*$, is normal if and only if
 $f$ induces a Galois extension  from the meromorphic function field $\mathcal{M}(\mathfrak{M}^*)$ to $\mathcal{M}(\mathfrak{M})$.
 Moreover, the Galois group is isomorphic to the covering transformation group in a natural manner.

 In \S\,154, Teichm\"uller defines the notion of normal covering between  principal regions. He then checks that for a normal covering from $\mathfrak{H}$ to $\mathfrak{H}^*$, with normal covering
 transformation group $\mathfrak{S}$, the regular quadratic differentials on $\mathfrak{H}^*$
 are exactly the regular quadratic differentials on $\mathfrak{H}$ that are invariant by $\mathfrak{S}$.
 
 The result in \S\,155 is a sequel to \S\,154:  The normal covering gives rise to a submanifold $R^*$ which is isometrically embedded as a geodesic
manifold and group-invariant.

 In \S\,156, Teichm\"uller applies the above discussion to study the action of conformal self-mappings  infinitesimally, that is, he works on the cotangent space
 of $R^\sigma$.

In \S\,157, he gives an example of a locally  extremal quasiconformal
mapping between tori which is not globally extremal. The existence of such a mapping
for  general surfaces
is conjectured in \S\,73.

In \S\,158, Teichm\"uller asks the following question:  is there a torus in
$3$-dimensional Euclidean space whose period ratio is not purely imaginary?\footnote{That is, the torus can be realized as $\mathbb{C}/<z\mapsto z+1, z\mapsto z+\tau>$, with $\mathrm{Im} \tau>0$.}
 He uses a heuristic argument to motivate an answer.

 \section{Generalization}

In \S\,159, Teichm\"uller considers generalizations of the questions on quasiconformal mappings raised in the previous sections to the case where the surface is the unit disc with the boundary fixed pointwise. He declares that one shall extend to a  quasiconformal mapping of the disc a ``sufficiently regular" topological mapping from the circle $\vert z\vert =1$ onto itself, and a few lines later, that ``certain regularity assumptions on the boudary must definitely be made whose appropriate formulation is a problem in itself."\footnote{By a result of Ahlfors and Beurling in \cite{ahlfors&beurling}, the regularity assumption is that the homeomorphism has to be a \textit{quasisymmetric map.}} The theory outlined in this paragraph is that of the ``non-reduced" Teichm\"uller space.\index{Teichm\"uller space!non-reduced}\index{non-reduced!Teichm\"uller space} He mentions that ``this is as if all boundary points were distinguished."  He discusses the role of quadratic differentials in this setting and the regularity properties at the boundary that are needed. He then makes the following assertion, which remained an open question for a long time:
\begin{quote}
Whenever a given boundary mapping can be extended to a quasiconformal mapping, then it can also be extended to an extremal quasiconformal mapping. It has the indicated form and 
$$
\iint_{\vert z \vert <1} \boxed{d\zeta}
$$
converges. If $\iint{\boxed{d\zeta}}$ converges, $d\zeta^2$ corresponds for each $K\geq 1$ to such a boundary mapping.
\end{quote} 
In the case of the unit disc, we can consider a holomorphic function instead of quadratic differential. We know now that there exist extremal mappings given by a constant $K$ and a quadratic differential $d\zeta^2$ satisfying
 $$
 \iint_{\vert z \vert <1} \boxed{d\zeta} =+\infty.
 $$ 
 Moreover, there also exist extremal mappings which are not given by a holomorphic function and whose associated Beltrami coefficient is not of constant absolute value. For these results, we refer the reader to \cite{yougos} (see also \cite{reichex}) and also to \cite{reichsurvey} for a global survey on this theory. 
 
 Teichm\"uller says that his aim is to stimulate further research on these questions. 
 
 He also considers arbitrary bordered  principal regions  with distinguished arcs on the boundary, that is, such that every point of such an arc is considered as a distinguished point. Let us note by the way that his later paper \cite{T31} deals with a problem that belongs to this setting.  More precisely, he considers there the problem of minimazing the quasiconformal dilatation for quasiconformal maps from the unit disc to itself which keep the boundary pointwise fixed and send $0$ to $-x$ ($0< x<1$).  Cf. also the commentary \cite{T31C} and Section 5.5.5 of \cite{reichsurvey}. 
 
 \begin{remark}  For a long time, the general form of extremal maps in the non-reduced theory was unknown. Indeed, it follows from work of Strebel that the Teichm\"uller extremal quasiconformal mappings\index{Teichm\"uller extremal mapping}\index{extremal mapping!Teichm\"uller} $f$ are also extremal in this setting. We recall here that these are quasiconformal maps whose associated Beltrami form $\mu_f$ is of the form
 $$
 k\cdot \frac{\overline{\varphi}}{\vert \varphi \vert},
 $$
where $0<k<1$ and $\varphi : \mathbb{D} \rightarrow \mathbb{D}$ is holomorphic and integrable on $\mathbb{D}$. As Reich writes in his survey \cite{reichsurvey}, Teichm\"uller suggested that \textbf{all} the extremal maps, for a certain boundary condition, are of this type. It was Strebel who constructed an extremal map\index{extremal mapping!non-Teichm\"uller} which is ``almost" of a 
  Teichm\"uller type, that is, its Beltrami form is of the given type except that the associated holomorphic function is not integrable.  A major advancement was made by V. Bo\v{z}in, N. Laki\'c, V. Markovi\'c and M. Mateljevi\'c who showed in \cite{yougos} that there exist extremal mappings (and even, uniquely extremal) which are not of Teichm\"uller's type, of finite or infinite norm.  
   \end{remark}
 \section{A metric}

This section contains only \S\,160.  As in the previous section, Teichm\"uller continues to explain some problems. He says that if we take a principal region $\mathfrak{h}$ with $\rho=0$ (negative Euler characteristic) and a variable point $\mathfrak{p}$ on $\mathfrak{h}$, then the set of pairs  $(\mathfrak{h},\mathfrak{p})$ defines a two-dimensional manifold.  He also declares that  there exists a Finsler metric on this manifold, where the distance between two points $(\mathfrak{h},\mathfrak{p})$ and $(\mathfrak{h}',\mathfrak{p}')$ is
defined by considering all the quasiconformal self-mappings of $\mathfrak{h}$ sending $\mathfrak{p}$ to $\mathfrak{p}'$, and taking the logarithm of the infimum of the quasiconformal dilatations of these maps. We know now that this question is related to the theory of the \textit{Bers fiber space}\index{Bers fiber space}\index{fiber space!Bers} and the \textit{forgetful map}. Indeed, in the case where the principal region $\mathfrak{h}$ is hyperbolic without boundary, the metric space that Teichm\"uller considers is the Bers fiber of $\mathfrak{h}$. Kra proved in \cite{kra} that if  the surface is not of exceptional type, then this two-dimensionnal manifold is not a Teichm\"uller disc (this result was also proved by Nag in \cite{nag}). Moreover, Kra proved that the Teichm\"uller distance on the fiber is not the hyperbolic one. Nevertheless, as Teichm\"uller recalls, this problem in the case of the sphere with three distinguished points was already solved in \S\,27, where the metric obtained is the hyperbolic metric. In fact, using the forgetful map, this fact is trivial because the fiber is the Teichm\"uller space of the sphere minus four points, which is exactly the hyperbolic disc.  Teichm\"uller   says that one can make a similar definition for distinguished points on the boundary. The question of studying this metric is presented as an open problem. The question is considered again in his paper  \cite{T31}, cf. also the commentary \cite{T31C}.

\section{An estimate}

In \S\,161, Teichm\"uller writes: ``We still have to deal with a question only researchers in the field should find some interest in. It concerns the proper development of the Gr\"otzsch-Ahlfors method."\index{length-area method}\index{Gr\"otzsch-Ahlfors method} He considers the problem of finding a coarse lower bound 
for the maximal dilatation of a  quasiconformal mapping between two given hexagons with distinguished pointing-out vertices. 
This lower bound is given in terms of side-lengths of the hexagons (see \S\,163),
by using a lemma provided in \S\,162. The lemma in \S\,162 is an improvement of the Gr\"otzsch-Ahlfors length-area\index{length-area method}\index{Gr\"otzsch-Ahlfors method} method. It gives a lower bound for the dilatation in terms 
of the average intersection number between the image of a horizontal strip with a vertical strip.
Note that a more precise lower bound is given in \S\,129.

\section{An extremal problem for conformal mappings}

This section starts at \S\,164. Teichm\"uller addresses the following question: \emph{Why do we study quasiconformal mappings?} He declares: ``Initially, quasiconformal mappings undoubtly appeared just as a generalization of conformal mappings."
He mentions the early work of Tissot\footnote{Nicolas Auguste Tissot\index{Tissot, Nicolas Auguste}  (1824-1897) is a famous cartographer, who also taught mathematics at the \'Ecole Polytechnique.} on the drawing of geographical maps, and works of Picard and Ahlfors. He says that these authors proved theorems on these mappings which were direct generalizations of theorems on conformal mappings. They were nevertheless interested in conformal mappings, and their results, he says, ``taught us to understand in what way the hypothesis of conformality restricts the behavior of a mapping." He adds (\S\,165): ``\emph{This view is outdated.} 
Today quasiconformal mappings are systematically used in studying conformal mappings." In particular, the fact of approximating a conformal mapping by a sequence of quasiconformal ones with dilatation quotients converging to 1 is not the essential part of the theory. Rather, it is important to construct, in the study of conformal mappings from a Riemann surface $\mathfrak{A}$ to a Riemann surface $\mathfrak{B}$, an explicit quasiconformal mapping from $\mathfrak{B}$ onto another Riemann surface  $\mathfrak{C}$ which is related to  $\mathfrak{B}$, and to study the triple  $(\mathfrak{A},  \mathfrak{B},  \mathfrak{C})$. He elaborates on this and he mentions applications to the so-called \emph{type problem}.\index{type problem}\index{problem!type} This expression designates the problem of determining whether a simply connected Riemann surface $\mathfrak{A}$ can be mapped holomorphically onto the plane or onto the unit disc, that is, the problem of determining whether the surface is ``parabolic" or ``hyperbolic." The idea is based on the fact that if a  quasiconformal mapping can be constructed from $\mathfrak{A}$ onto the plane (respectively the unit disc), then $\mathfrak{A}$ can be mapped holomorphically onto the plane (respectively the unit disc). Teichm\"uller mentions Picard's theorem\index{theorem!Picard} in this theory. This question is also related to Nevanlinna's theory.\footnote{Rolf Nevanlinna\index{Nevanlinna, Rolf}  (1895-1980) was one of Teichm\"uller's teachers at G\"ottingen. In the preface to his book 
\emph{Analytic functions} \cite{Nevanlinna}, Nevanlinna writes: ``The present monograph on analytic functions coincides to a large extent with the presentation of the modern theory of single-value analytic functions given in my earlier works ``Le th\'eor\`eme de Picard-Borel et la th\'eorie des fonctions m\'eromorphes" \cite{Nevanlinna1} and ``Eindeutige analytische Funktionen" \cite{Nevanlinna}. In Section XII, titled
 \emph{The type of a Riemann surface}, Nevanlinna mentions a question on the type problem\index{type problem}\index{problem!type} answered by Teichm\"uller in his paper \cite{T200} (p. 312 of the English edition).} 
 Teichm\"uller mentions works of other of his contemporary mathematicians.\footnote{Teichm\"uller refers to the paper ``Zum Umkehrproblem der Wertvertelungslehre" \cite{Ullrich}, by E. Ullrich and to the \emph{habilitationsschrift} 
 of H. Wittich. In the first work, the author studies a new class of Riemann surfaces called ``surfaces with finitely many periodic ends." He proves the existence of meromorphic functions which have prescribed rational valued \emph{deficiency} of \emph{ramification index} at finitely many points, and sum 2. In the second work, the author  uses quasi-conformal mappings in his study of value-distribution theory, that is, the theory (also called Nevanlinna's theory) whose aim is to count the number of times a holomorphic function $f(z)$ assumes a certain value as $z$ grows in size. This theory is a wide generalization of the classical Picard  theorem.}  Then he states that his ``extremal problem for quasiconformal mappings and its solution can lead to heuristic methods for solving extremal problems about conformal mappings whose nature was initially completely different." He then returns to some results he already discussed in the present paper, in particular the ``Gr\"otzsch-Ahlfors method,"\index{Gr\"otzsch-Ahlfors method}\index{length-area method} and he declares: ``Some connection between the extremal problems of conformal or quasiconformal mappings and quadratic differentials can also be noticed in the work of Gr\"otzsch."
He concludes his paper with the announcement of further results in future papers, in particular, detailed proofs of results announced in the present paper. 

The interested reader may find several other of Teichm\"uller's papers translated and commented on in this Handbook.

\printindex
\end{document}